\documentclass[12pt,twoside]{article} 
\usepackage{amsmath,amssymb,latexsym,theorem,bbm,shapepar,natbib,epsfig}
\newfont{\msa}{msam10 scaled\magstep1}
\newfont{\ssmsa}{msam9}

\def\qlim{\mathop{\hbox{\rm l.i.m.}}}

\newcommand{\proofend}{\hfill$\square$}

\numberwithin{equation}{section}

\theorembodyfont{\em}
\newtheorem{Lem}{Lemma}[section]
\newtheorem{Thm}[Lem]{Theorem}

\theorembodyfont{\rm}
\newtheorem{Rem}[Lem]{Remark}
\newtheorem{Ex}[Lem]{Example}
   
\title{Parameter estimation in linear regression driven by a
  Gaussian sheet} 

\author{{\sc S\'andor Baran} and {\sc Kinga Sikolya}\\
         Faculty of Informatics, University of Debrecen\\
         Kassai \'ut 26, H--4028 Debrecen, Hungary}

\date{}

\begin{document}

\pagestyle{myheadings}
\markboth{\rm Baran, Sikolya}
         {\rm Parameter estimation in linear regression driven by a
           Gaussian sheet} 

\maketitle

\begin{abstract}
The problem of estimating the parameters of a linear
regression model \ $Z(s,t)=m_1g_1(s,t)+ \cdots +
m_pg_p(s,t)+U(s,t)$ \  based on observations of \ $Z$ \ on a spatial
domain \ $G$ \ of special shape is considered,  where the driving
process \ $U$ \ is a 
Gaussian random field and \ $g_1, \ldots, g_p$ \ are known
functions. Explicit forms of the maximum-likelihood estimators of
the parameters are derived in the cases when \ $U$ \ is either a Wiener or a
stationary or nonstationary Ornstein-Uhlenbeck sheet. Simulation
results are also presented,
where the driving random sheets are simulated with the help of their
Karhunen-Lo\`eve expansions.


\end{abstract}

\section{Introduction}
   \label{sec:sec1}

The Wiener sheet is one of the most important examples of Gaussian
random fields. It has various applications in statistical
modeling. Wiener sheet appears as limiting process of some
random fields defined on the interface of the Ising model \citep{kuma},
it is used to model random polymers \citep{doug}, to
describe the dynamics of Heath--Jarrow--Morton type forward interest rate models
\citep{gold} or to model random mortality surfaces
\citep{bimi}. Further, \citet{carter} considers the problem of
estimation of the mean in a nonparametric regression on a two-dimensional
regular grid of design points and constructs a Wiener sheet process on
the unit square with a drift that is almost the mean function in the
nonparametric regression. 

The stationary Ornstein-Uhlenbeck sheet \
$\{\widetilde X (s,t):s,t\in{\mathbb R}\}$ 
 \ is a zero mean Gaussian process with covariance structure
\begin{equation}
   \label{ou1}
 {\mathsf E}\widetilde{X}(s_1,t_1)\widetilde{X}(s_2,t_2)
   =\frac{\sigma^2}{4\alpha\beta}{\mathrm
     e}^{-\alpha|s_2-s_1|-\beta|t_2-t_1|},
\end{equation}
 where \ $\alpha>0$, \ $\beta>0$, \ $\sigma>0$, \ while the random
 field
\begin{equation}
    \label{ou2}
 X(s,t)=\sigma\int\limits_0^s\int\limits_0^t{\mathrm
   e}^{\alpha(u-s)+\beta(v-t)}\,{\mathrm d} W(u,v), 
   \qquad s,t\geq0,
\end{equation}
 where \ $\alpha\in{\mathbb R}$, \ $\beta\in{\mathbb R}$, \ $\sigma>0$ \ and
 \ $\{W(s,t):s,t\geq0\}$ \ is a standard Wiener sheet, can be considered as the
 Ornstein-Uhlenbeck sheet with zero initial condition on the axes.

Ornstein-Uhlenbeck sheets play 
role e.g. in potential
theory \citep{fepra} and,  similarly to the Wiener sheet, they also appear as
driving fields in forward interest rate models \citep{gold,scs}.

In this paper we consider a linear regression model driven by a 
Gaussian sheet, that is a random field  \ 
\begin{equation}
   \label{model}
Z(s,t):=m_1g_1(s,t)+ \cdots +
m_pg_p(s,t)+U(s,t)
\end{equation}
observed on a domain \ $G$, \ where \ $g_1, \ldots, g_p$ \ are known
functions and \ $U$ \ is either a Wiener or an
Ornstein-Uhlenbeck sheet, and we
determine the maximum likelihood estimator (MLE) of the unknown
parameters \ $m_1,\ldots, m_p$. 

In principle, the Radon-Nikodym derivative of Gaussian measures might be
derived from the general Feldman-Hajek theorem \citep{kuo}, but in
most of the cases explicit calculations can not be carried out. In the case
when \ $U$ \ is a standard Wiener
sheet, $p=1$ \ and \ $g_1\equiv
1$ \ (shifted Wiener sheet), \  the MLE of the unknown parameter is
given e.g. in 
\citet{rosa1}, where the estimator is expressed as a function of a
usually unknown random variable satisfying some characterizing equation.
In several cases the exact form of this random variable can be derived
by a method proposed by \citet{rosa2}, based on linear stochastic
partial differential equations. \citet{aratonm1} used Rozanov's method
to find the MLE of the shift parameter of a shifted Wiener sheet
observed on a special domain. \citet{bpz1} considered the model of
\citet{aratonm1}, and applying an essentially simpler direct discrete
approximation approach the authors found the MLE of the shift parameter under
much weaker conditions. Later this discrete approximation was used to
determine the MLE of the unknown parameter for the model \eqref{model}
with \ $p=1$ \ and a more complicated domain of 
observations, when \ $U$ \ is a standard Wiener sheet.

\citet{aratonm2} also studied the case when \ $U$ \ is an
Ornstein-Uhlenbeck sheet, \ $p=1$ \ and \ $g_1\equiv 1$, \ and using
partial stochastic differential equations found
the MLE of the unknown parameter based on the observation of the
random field \ $Z$ \ on a rectangular domain. This result was
generalized by \citet{bpz2} for the case \ $p=1$ \ and \ $g_1$ \ with
slight analytic restrictions.

In the present paper we consider the same type of domain
\ $G$ \ as in \citet{bpz3} and extend their result for the general
model \eqref{model}. We also consider the cases when the driving
process \ $U$ \ is a stationary and a zero start Ornstein-Uhlenbeck
sheet and generalize the results of \citet{aratonm2} and
\citet{bpz3}. Moreover, 
we present some simulation results to illustrate the theoretical ones
where the driving Gaussian random sheets are simulated with the help
of their  Karhunen-Lo\`eve expansions. 
The proofs of the theorems are given in the Appendix.

\section{Models and estimators}
   \label{sec:sec2}

Consider the model \eqref{model} with some given functions \  
$g_1, \ldots, g_p:{\mathbb R}^2_+ \to {\mathbb R}$ \ and with 
unknown regression parameters \ 
$m_1,\ldots ,m_p\in{\mathbb R}$. \ Let \ $[a,c]\subset (0,\infty)$ \ and \
$b_1,b_2\in (a,c)$, \ let \ $\gamma _{1,2}:[a,b_1]\to {\mathbb R}$ \
and \ $\gamma _0:[b_2,c]\to {\mathbb R}$ \ be continuous, strictly
decreasing functions and let \ $\gamma _1:[b_1,c]\to {\mathbb R}$ \
and \ $\gamma _2:[a,b_2]\to {\mathbb R}$ \ be continuous, strictly
increasing functions with \ $\gamma _{1,2}(b_1)=\gamma _1(b_1)>0$, \ 
$\gamma _2(b_2)=\gamma _0(b_2)$, \ $\gamma _{1,2}(a)=\gamma _2(a)$ \
and \ $\gamma _1(c)=\gamma _0(c)$. \ \ Consider the curve
 \ $\Gamma:=\Gamma_{1,2}\cup\Gamma_1\cup\Gamma_2\cup\Gamma_0$,
 \ where
 \begin{alignat*}{2}
  \Gamma _{1,2}&:=\Big\{\big(s,\gamma_{1,2}(s)\big):s\in [a,b_1]\Big\}, \qquad
  &\Gamma _1:=\Big\{\big(s,\gamma_1(s)\big):s\in [b_1,c]\Big\}, \\
  \Gamma _2&:=\Big\{\big(s,\gamma_2(s)\big):s\in [a,b_2]\Big\}, 
  &\Gamma _0:=\Big\{\big(s,\gamma_0(s)\big):s\in [b_2,c]\Big\},
 \end{alignat*}
 and for a given \ $\varepsilon>0$ \ let \ $\Gamma_{1,2}^\varepsilon$,
 \ $\Gamma_1^\varepsilon$, 
 \ $\Gamma_2^\varepsilon$ \ and \ $\Gamma_0^\varepsilon$ \ denote the
 \ inner $\varepsilon$-strip 
 of \ $\Gamma_{1,2}$, \ $\Gamma_1$, \ $\Gamma_2$ \ and \ $\Gamma_0$,
 \ respectively, that is e.g.
 \begin{align*}
  \Gamma_{1,2}^{\varepsilon} :=\big\{(s,t) \in {\mathbb R}^2
      :\, &s\in [a,a+\varepsilon], \ t\in
      [\gamma_{1,2}(s),\gamma_{1,2}(a)] \ \text{or} \\ 
        &s\in [a+\varepsilon,b_1], \
           t\in[\gamma_{1,2}(s),\gamma_{1,2}(s)+\varepsilon] \big\}.   
 \end{align*}
Suppose that there exists an \ $\varepsilon>0$ \ such that
 \begin{equation}
    \label{eq:eq2.1}
  \Gamma_1^\varepsilon\cap\Gamma_2^\varepsilon=\emptyset \quad \text{and} \quad
  \Gamma_{1,2}^\varepsilon\cap\Gamma_0^\varepsilon=\emptyset,
 \end{equation}
 and consider the set \ $G:=G_1\cup G_2\cup G_3$, \ where
 \begin{align*}
  G_1&:=\big\{(s,t)\in{\mathbb R}^2
              :s\in[a,b_1\land b_2],\,
               t\in[\gamma_{1,2}(s),\gamma_2(s)]\big\},\\
  G_2&:=\begin{cases}
         \big\{(s,t)\in{\mathbb R}^2
               :s\in[b_1,b_2],\,t\in[\gamma_1(s),\gamma_2(s)]\big\},
          &\quad\text{if \ $b_1\leq b_2$,}\\
         \big\{(s,t)\in{\mathbb R}^2
               :s\in [b_2,b_1],\,t\in[\gamma_{1,2}(s),\gamma_0(s)]\big\},
          &\quad\text{if \ $b_1> b_2$,}
        \end{cases}\\
  G_3&:=\big\{(s,t)\in{\mathbb R}^2
              :s\in[b_1\lor b_2,c],\,t\in[\gamma_1(s),\gamma_0(s)]\big\}. 
 \end{align*}
\begin{figure}[tb]
 \begin{center}
\leavevmode
\epsfig{file=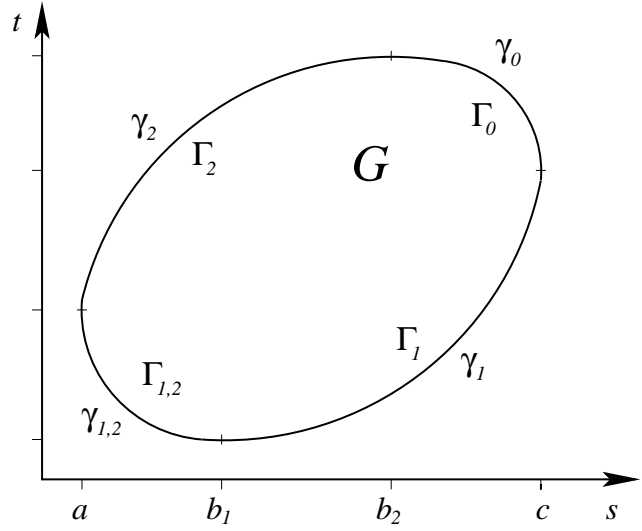,height=7cm}
\caption{\label{fig:fig1} An example of a set of observations \ $G$.}  
  \end{center}  
\end{figure}
An example of such a set of observations can be seen of Figure
\ref{fig:fig1}.

First we study the case when \ $U$ \ a standard Wiener sheet and consider
the random field \ $Z(s,t):=m_1g_1(s,t)+ \cdots +
m_pg_p(s,t)+W(s,t)$. \ The following theorem is an obvious extension
of Theorem 2.1 of \citet{bpz3} and can be proved exactly in the same
way. The proof is based on the discrete approximation method described
in  \citet{bpz2,bpz1, bpz3}, which relies on the results
of \citet[Section 2.3.2]{arato}. 

\begin{Thm}
   \label{mainW} \ 
If \ $g_1, \ldots, g_p$ \ are twice continuously differentiable 
inside \  $G$ \ and the partial derivatives  \ $\partial_1g_i, \ 
\partial_2g_i$ \ and \ $\partial_1\partial_2g_i, \ i=1, \ldots, p$, \ 
 can be continuously extended to  $G$  then  the probability measures
 \ ${\mathsf P}_Z$ \ and \ ${\mathsf P}_W$, \ generated on \ $C(G)$ \
 by the sheets \ $Z$ \ and 
 \ $W$, \ respectively, are equivalent and the Radon-Nikodym derivative of 
 \ ${\mathsf P}_Z$ \ with respect to \ ${\mathsf P}_W$ \ equals
 \begin{equation*}
  \frac{{\mathrm d}{\mathsf P}_Z}{{\mathrm d}{\mathsf
      P}_W}(Z)=\exp\left\{-\frac 12\big({\mathbf m}^{\top}A{\mathbf
      m}-2\zeta ^{\top} {\mathbf m}\big)\right\}, 
 \end{equation*}
where \ $A:=\big (A_{k,\ell} \big)_{k,\ell=1}^p$,  \ ${\mathbf
  m}:=(m_1,\ldots, m_p)^{\top}$ \ and \ 
  $\zeta :=\big(\zeta_1, \ldots , \zeta_p\big)^{\top}$ \ with
 \begin{align}
   \label{eq:eq2.2}
  A_{k,\ell}:=&\frac{g_k\big(b_1,\gamma_{1,2}(b_1)\big)\,
    g_{\ell}\big(b_1,\gamma_{1,2}(b_1)\big)}{b_1\gamma_{1,2}(b_1)} \\  
      &+\int\limits_a^{b_1}
        \frac{\big[g_k\big(s,\gamma_{1,2}(s)\big)-
          s\partial_1g_k\big(s,\gamma_{1,2}(s)\big)\big] 
       \big[g_{\ell}\big(s,\gamma_{1,2}(s)\big)-
       s\partial_1g_{\ell}\big(s,\gamma_{1,2}(s)\big)\big]} 
             {s^2\gamma_{1,2}(s)}
        \,{\mathrm d} s \nonumber \\
      &+\int\limits_{b_1}^c
        \frac{\partial_1 g_k\big(s,\gamma_1(s)\big)\,\partial_1
          g_{\ell}\big(s,\gamma_1(s)\big)} 
             {\gamma_1(s)}
        \,{\mathrm d} s
     +\int\limits_{\gamma_2(a)}^{\gamma_2(b_2)}
        \frac{\partial_2g_k\big(\gamma_2^{-1}(t),t\big)\,
          \partial_2g_{\ell}\big(\gamma_2^{-1}(t),t\big)} 
             {\gamma_2^{-1}(t)}
        \,{\mathrm d} t \nonumber \\
     & +\int\limits_{\gamma_{1,2}(b_1)}^{\gamma_{1,2}(a)}
        \frac{\partial_2g_k\big(\gamma_{1,2}^{-1}(t),t\big)\,
          \partial_2g_{\ell}\big(\gamma_{1,2}^{-1}(t),t\big) }
             {\gamma_{1,2}^{-1}(t)}
        \,{\mathrm d} t
      +\iint\limits_G
        \partial_1\partial_2g_k(s,t)\,
        \partial_1\partial_2g_{\ell}(s,t){\mathrm d}
        s\,{\mathrm d} t,  \nonumber
 \end{align}
 and
 \begin{align}
    \label{eq:eq2.3}
  \zeta_k:=&\frac{g_k\big(b_1,\gamma_{1,2}(b_1)\big)
    Z\big(b_1,\gamma_{1,2}(b_1)\big)}   
               {b_1\gamma_{1,2}(b_1)}
         +\int\limits_{b_1}^c
           \frac{\partial_1 g_k\big(s,\gamma_1(s)\big)}
                {\gamma_1(s)}
           \,Z\big({\mathrm d} s,\gamma_1(s)\big) \\         
         &+\int\limits_a^{b_1}
            \frac{\big[g_k\big(s,\gamma_{1,2}(s)\big)
              -s\partial_1g_k\big(s,\gamma_{1,2}(s)\big)\big]} 
                 {s^2\gamma_{1,2}(s)} 
            \Big[Z\big(s,\gamma_{1,2}(s)\big)\,{\mathrm d} s-sZ\big({\mathrm d}
            s,\gamma_{1,2}(s)\big)\Big] \nonumber \\ 
         &+\int\limits_{\gamma_2(a)}^{\gamma_2(b_2)}
            \frac{\partial_2g_k\big(\gamma_2^{-1}(t),t\big)}
                 {\gamma_2^{-1}(t)}
            \,Z\big(\gamma_2^{-1}(t),{\mathrm d} t\big)
          +\int\limits_{\gamma_{1,2}(b_1)}^{\gamma_{1,2}(a)}
            \frac{\partial_2g_k\big(\gamma_{1,2}^{-1}(t),t\big)}
                 {\gamma_{1,2}^{-1}(t)}
            \,Z\big(\gamma_{1,2}^{-1}(t),{\mathrm d} t\big)\nonumber \\
         &+\iint\limits_G
            \partial_1\partial_2g_k(s,t)\,Z({\mathrm d} s,{\mathrm d}
            t). \nonumber 
 \end{align}
The maximum likelihood estimator of the parameter vector \ ${\mathbf
  m}$ \ based on the 
 observations \ $\{Z(s,t):(s,t)\in G\}$ \ has the form \ ${\widehat
   {\mathbf m}}=A^{-1}\zeta$ \ and 
 has a $p$-dimensional normal distribution with mean \ ${\mathbf m}$ \ and
 covariance matrix \ $A^{-1}$. 
\end{Thm}

\begin{Rem}
  \label{integral} \
We remark that the weighted $L^2$-Riemann integrals of partial
derivatives of the Wiener sheet (and of other  $L^2$-processes) along a
curve are defined in the sense of \citet[Definition 4.1]{bpz3}. This
means that if \ $Z$ \ is an $L^2$-process given along an
$\varepsilon$-neighborhood of a curve 
 \ $\Gamma:=\big\{(s,\gamma(s)):s\in[a,b]\big\}$, \ where
 \ $\gamma:[a,b]\to{\mathbb R}$ \ is strictly monotone and \
 $y:[a,b]\to{\mathbb R}$ \ is a function, then
 \begin{align*}
  \int\limits_a^by(s)\,Z({\mathrm d} s,\gamma(s))
  &:=\qlim_{h\to 0}
      \frac 1h
      \int\limits_a^by(s)\big[Z(s+h,\gamma(s))-Z(s,\gamma(s))\big]\,
      {\mathrm d} s, \\     
  \int\limits_{\gamma(a)}^{\gamma(b)}
   y(\gamma^{-1}(t))\,Z(\gamma^{-1}(t),{\mathrm d} t)
  &:=\qlim_{h\to 0} 
      \frac 1h 
      \int\limits_{\gamma(a)}^{\gamma(b)}
       y(\gamma^{-1}(t))\big[Z(\gamma^{-1}(t),t+h)-
   Z(\gamma^{-1}(t),t)\big]\,{\mathrm d} t,  
 \end{align*}
 if the right hand sides exist. 
\end{Rem}

 The stationary Ornstein-Uhlenbeck sheet \ $\{\widetilde X(s,t):s,t\in
{\mathbb R}\}$ \ defined by covariance structure \eqref{ou1} can
be represented as 
 \begin{equation} 
      \label{eq:eq2.4}
\widetilde X(s,t)
   =\frac{\sigma}{2\sqrt{\alpha\beta}}\,{\mathrm e}^{-\alpha s-\beta t}\,
    W\big({\mathrm e}^{2\alpha s},{\mathrm e}^{2\beta t}\big),
   \qquad s,t\in{\mathbb R}.
 \end{equation}
Consider the sheet \ $Z(s,t):=m_1h_1(s,t)+ \cdots +
m_ph_p(s,t)+\widetilde X(s,t), \ (s,t)\in {\mathbb R}_+$. \
\ Applying Theorem \ref{mainW} for the functions
 \begin{equation}
    \label{eq:eq2.5}
  g_k(u,v)
  =\frac{2\sqrt{\alpha\beta uv}}{\sigma}\,
   h_k\left(\frac{\log{u}}{2\alpha},\frac{\log{v}}{2\beta}\right),
   \quad k=1,2,\ldots ,p,
 \end{equation}
and for the domain \ $\widetilde G$ \ bounded by the curve
 \ $\widetilde\Gamma:=\widetilde\Gamma_{1,2}\cup\widetilde\Gamma_1
 \cup\widetilde\Gamma_2\cup\widetilde\Gamma_0$, 
 \ where
 \begin{alignat*}{2}
  \widetilde\Gamma _{1,2}&:=\Big\{\big(u,\widetilde\gamma_{1,2}(u)\big):u\in
  \big[{\mathrm e}^{2\alpha a},{\mathrm e}^{2\alpha b_1}\big]\Big\},
  \qquad
&\widetilde\Gamma _1:=\Big\{\big(u,\widetilde\gamma_1(u)\big):u\in
  \big[{\mathrm e}^{2\alpha b_1},{\mathrm e}^{2\alpha c}\big]\Big\},
  \\
 \widetilde\Gamma _2&:=\Big\{\big(u,\widetilde\gamma_2(u)\big):u\in
  \big[{\mathrm e}^{2\alpha a},{\mathrm e}^{2\alpha b_2}\big]\Big\},
  \qquad
&\widetilde\Gamma _0:=\Big\{\big(u,\widetilde\gamma_0(u)\big):u\in
  \big[{\mathrm e}^{2\alpha b_2},{\mathrm e}^{2\alpha c}\big]\Big\}, 
 \end{alignat*}
with
\begin{equation*}
\widetilde\gamma_j(u):=\exp \Big(2\beta \gamma_j\big(\log
(u)/2\alpha\big)\Big), \qquad j\in \big\{\{1,2\},\{1\},\{2\},\{0\}\big\},
\end{equation*}
we obtain the following result.

\begin{Thm}
   \label{mainOU1} \
If \ $h_1, \ldots, h_p$ \ are twice continuously differentiable 
inside \  $G$ \ and the partial derivatives  \ $\partial_1h_i, \ 
\partial_2h_i$ \ and \ $\partial_1\partial_2h_i, \ i=1, \ldots, p$, \ 
 can be continuously extended to  $G$  then the probability measures
 \ ${\mathsf P}_Z$ \ and \ ${\mathsf P}_{\widetilde X}$, \ generated
 on \ $C(G)$ \ by the sheets \ $Z$ \ and 
 \ $\widetilde X$, \ respectively, are equivalent and the
 Radon-Nikodym derivative of  
 \ ${\mathsf P}_Z$ \ with respect to \ ${\mathsf P}_{\widetilde X}$ \ equals
 \begin{equation*}
  \frac{{\mathrm d}{\mathsf P}_Z}{{\mathrm d}{\mathsf
      P}_{\widetilde X}}(Z)=\exp\left\{-\frac {\alpha\beta}{2\sigma
      ^2}\big({\mathbf m}^{\top}A{\mathbf 
      m}-2\zeta ^{\top} {\mathbf m}\big)\right\}, 
 \end{equation*}
where \ $A:=\big (A_{k,\ell} \big)_{k,\ell=1}^p$,  \ ${\mathbf
  m}:=(m_1,\ldots, m_p)^{\top}$ \ and \ 
  $\zeta :=\big(\zeta_1, \ldots , \zeta_p\big)^{\top}$ \ with
\begin{align}
  \label{eq:eq2.6}
A_{k,\ell}:=&\,h_k\big(a,\gamma_2(a)\big)h_{\ell}\big(a,\gamma_2(a)\big)
+ h_k\big(c,\gamma_1(c)\big)h_{\ell}\big(c,\gamma_1(c)\big) \\
& + h_k\big(b_1,\gamma_1(b_1)\big)h_{\ell}\big(b_1,\gamma_1(b_1)\big)
 + h_k\big(b_2,\gamma_2(b_2)\big)h_{\ell}\big(b_2,\gamma_2(b_2)\big)
 \nonumber\\
 &+\int\limits_a^{b_1}\Big[\alpha
 h_k\big(s,\gamma_{1,2}(s)\big)h_{\ell}\big(s,\gamma_{1,2}(s)\big)
 +\alpha^{-1}\partial_1h_k\big(s,\gamma_{1,2}(s)\big) 
\partial_1h_{\ell}\big(s,\gamma_{1,2}(s)\big)\Big]\,{\mathrm d}s \nonumber\\
&+\int\limits_{b_1}^c\Big[\alpha h_k\big(s,\gamma_1(s)\big)+\partial
_1h_k\big(s,\gamma_1(s)\big)\Big]\Big[h_{\ell}\big(s,\gamma_1(s)\big) 
 +\alpha^{-1} \partial_1h_{\ell}\big(s,\gamma_1(s)\big)\Big]\,{\mathrm
   d}s \nonumber\\ 
&+\int\limits_a^{b_2}\Big[\alpha h_k\big(s,\gamma_2(s)\big)-\partial
_1h_k\big(s,\gamma_2(s)\big)\Big]\Big[h_{\ell}\big(s,\gamma_2(s)\big) 
 -\alpha^{-1} \partial_1h_{\ell}\big(s,\gamma_2(s)\big)\Big]\,{\mathrm
   d}s \nonumber\\  
&+\int\limits_{b_2}^c\Big[\alpha h_k\big(s,\gamma_0(s)\big)
h_{\ell}\big(s,\gamma_0(s)\big)+
\alpha^{-1}\partial_1h_k\big(s,\gamma_0(s)\big)
\partial_1h_{\ell}\big(s,\gamma_0(s)\big)\Big]\,{\mathrm d}s \nonumber
\\
&+\int\limits_{\gamma_{1,2}(b_1)}^{\gamma_{1,2}(a)}\Big[\beta 
 h_k\big(\gamma_{1,2}^{-1}(t),t\big)-\partial_2
 h_k\big(\gamma_{1,2}^{-1}(t),t\big)\Big]\Big[ 
 h_{\ell}\big(\gamma_{1,2}^{-1}(t),t\big)-\beta^{-1}
\partial_2 h_{\ell}\big(\gamma_{1,2}^{-1}(t),t\big)\Big]\,{\mathrm d}t
\nonumber \\
 &+\int\limits_{\gamma_1(b_1)}^{\gamma_1(c)}\Big[\beta
 h_k\big(\gamma_1^{-1}(t),t\big)h_{\ell}\big(\gamma_1^{-1}(t),t\big)
 +\beta^{-1}\partial_2h_k\big(\gamma_1^{-1}(t),t\big)
  \partial_2h_{\ell}\big(\gamma_1^{-1}(t),t\big)\!\Big]\,{\mathrm d}t 
\nonumber 
\end{align}
\begin{align}
\phantom {A_{k,\ell}:=}
 &+\int\limits_{\gamma_2(a)}^{\gamma_2(b_2)}\Big[\beta
 h_k\big(\gamma_2^{-1}(t),t\big)h_{\ell}\big(\gamma_2^{-1}(t),t\big)
 +\beta^{-1} \partial_2h_k\big(\gamma_2^{-1}(t),t\big) 
\partial_2h_{\ell}\big(\gamma_2^{-1}(t),t\big)\Big]\,{\mathrm d}t
\nonumber\\
&+\int\limits_{\gamma_0(c)}^{\gamma_0(b_2)} \Big[\beta 
 h_k\big(\gamma_0^{-1}(t),t\big)+\partial_2h_k\big(\gamma_0^{-1}(t),t\big)\Big]
 \Big[ h_{\ell}\big(\gamma_0^{-1}(t),t\big)+
\beta^{-1}\partial_2h_{\ell}\big(\gamma_0^{-1}(t),t\big)\Big]\,{\mathrm
  d}t \nonumber\\ 
&+\iint\limits_G \big [\alpha\beta
h_k(s,t)h_{\ell}(s,t)+\alpha^{-1}\beta \partial_1h_k(s,t)\partial_1h_{\ell}(s,t)
+\alpha\beta^{-1} \partial_2h_k(s,t)\partial_2h_{\ell}(s,t)\nonumber \\
&\phantom{+\iint\limits_G=}
+\alpha^{-1}\beta^{-1} \partial_1\partial_2h_k(s,t)
\partial_1\partial_2h_{\ell}(s,t) \big]\,{\mathrm d}s\, {\mathrm d}t  \nonumber
\end{align}
and
\begin{align}
\label{eq:eq2.7}
\zeta_k:=&\,4h_k\big(b_1,\gamma_1(b_1)\big)Z\big(b_1,\gamma_1(b_1)\big) \\
&+2\int\limits_{b_1}^c\Big[\alpha h_k\big(s,\gamma_1(s)\big)
 +\partial_1h_k\big(s,\gamma_1(s)\big)\Big]
\Big[Z\big(s,\gamma_1(s)\big){\mathrm d}s+\alpha^{-1}Z\big({\mathrm
  d}s,\gamma_1(s)\big)\Big]\nonumber\\
&+2\int\limits_a^{b_1}\Big[\alpha h_k\big(s,\gamma_{1,2}(s)\big)
 -\partial_1h_k\big(s,\gamma_{1,2}(s)\big)\Big]
\Big[Z\big(s,\gamma_{1,2}(s)\big){\mathrm d}s-\alpha^{-1}Z\big({\mathrm
  d}s,\gamma_{1,2}(s)\big)\Big]\nonumber\\
&+2\int\limits_{\gamma_2(a)}^{\gamma_2(b_2)}\Big[\beta
 h_k\big(\gamma_2^{-1}(t),t\big)+ \partial_2h_k\big(\gamma_2^{-1}(t),t\big)\Big]
\Big[Z\big(\gamma_2^{-1}(t),t\big){\mathrm d}t+
\beta^{-1}Z\big(\gamma_2^{-1}(t),{\mathrm d}t\big)\Big]\nonumber \\ 
&+2\int\limits_{\gamma_{1,2}(b_1)}^{\gamma_{1,2}(a)}\Big[\beta
 h_k\big(\gamma_{1,2}^{-1}(t),t\big)+ \partial_2h_k
 \big(\gamma_{1,2}^{-1}(t),t\big)\Big]  
\Big[Z\big(\gamma_{1,2}^{-1}(t),t\big){\mathrm d}t+
\beta^{-1}Z\big(\gamma_{1,2}^{-1}(t),{\mathrm d}t\big)\Big] \nonumber\\ 
&+\iint\limits_G \big [\alpha\beta
h_k(s,t)+\beta \partial_1h_k(s,t)+\alpha\partial_2h_k(s,t)+
\partial_1\partial_2h_k(s,t)\big] \nonumber\\
&\phantom{+\iint\limits_G=}\times \big [Z(s,t){\mathrm d}s\, {\mathrm d}t
+\alpha^{-1} Z({\mathrm d}s,t){\mathrm d}t  
+\beta^{-1} Z(s,{\mathrm d}t) {\mathrm d}s 
+\alpha^{-1}\beta^{-1} Z({\mathrm d}s,{\mathrm d}t)\big].  \nonumber
\end{align}
The maximum likelihood estimator of the parameter vector \ ${\mathbf
  m}$ \ based on the 
 observations \ $\{Z(s,t):(s,t)\in G\}$ \ has the form \ ${\widehat
   {\mathbf m}}=A^{-1}\zeta$ \ and 
 has a $p$-dimensional normal distribution with mean \ ${\mathbf m}$ \ and
 covariance matrix \ $A^{-1}$. 
\end{Thm}

Finally, consider the sheet \ $Z(s,t):=m_1h_1(s,t)+ \cdots +
m_ph_p(s,t)+X(s,t), \ (s,t)\in {\mathbb R}_+$, where 
 \ $\{X(s,t):s,t\geq0\}$ \ is the zero start  Ornstein--Uhlenbeck
 sheet defined by \eqref{ou2}. For  
\ $\alpha\not=0$ \ and \ $\beta\not=0$ \ the sheet \
$\{X(s,t):s,t\geq0\}$ \ can be characterized as 
 a zero mean Gaussian process with
\begin{equation*}
 {\mathsf E} X(s_1,t_1)X(s_2,t_2)
   =\frac{\sigma^2}{4\alpha\beta}
    \left({\mathrm e}^{-\alpha|s_1-s_2|}-{\mathrm e}^{-\alpha(s_1+s_2)}\right)
    \left({\mathrm e}^{-\beta|t_1-t_2|}-{\mathrm
        e}^{-\beta(t_1+t_2)}\right),
 \end{equation*}
 hence, for example, in case \ $\alpha>0$ \ and \ $\beta>0$ \ it can
 also be represented as
 \begin{equation}
    \label{eq:eq2.8}
  X(s,t)=\frac{\sigma}{2\sqrt{\alpha\beta}}\,{\mathrm e}^{-\alpha s-\beta t}\,
          W({\mathrm e}^{2\alpha s}-1,{\mathrm e}^{2\beta t}-1),
   \qquad s,t\geq0.
 \end{equation}
In this way similarly to the stationary case one can apply Theorem
\ref{mainW} for the functions
 \begin{equation}
    \label{eq:eq2.9}
  g_k(u,v)
  =\frac{2\sqrt{\alpha\beta (u+1)(v+1)}}{\sigma}\,
   h_k\left(\frac{\log (u+1)}{2\alpha},\frac{\log (v+1) }{2\beta}\right),
   \quad k=1,2,\ldots ,p,
 \end{equation}
and for the domain \ $\widehat G$ \ bounded by the curve
 \ $\widehat\Gamma:=\widehat\Gamma_{1,2}\cup\widehat\Gamma_1
 \cup\widehat\Gamma_2\cup\widehat\Gamma_0$, 
 \ where
 \begin{alignat*}{2}
  \widehat\Gamma _{1,2}&:=\!\Big\{\big(u,\widehat\gamma_{1,2}(u)\big):u\in
  \big[{\mathrm e}^{2\alpha a}\!-\!1,{\mathrm e}^{2\alpha b_1}\!-\!1\big]\Big\},
  \quad
&\widehat\Gamma _1:=\!\Big\{\big(u,\widehat\gamma_1(u)\big):u\in
  \big[{\mathrm e}^{2\alpha b_1}\!-\!1,{\mathrm e}^{2\alpha c}\!-\!1\big]\Big\},
  \\
 \widehat\Gamma _2&:=\!\Big\{\big(u,\widehat\gamma_2(u)\big):u\in
  \big[{\mathrm e}^{2\alpha a}\!-\!1,{\mathrm e}^{2\alpha b_2}\!-\!1\big]\Big\},
  \quad
&\widehat\Gamma _0:=\!\Big\{\big(u,\widehat\gamma_0(u)\big):u\in
  \big[{\mathrm e}^{2\alpha b_2}\!-\!1,{\mathrm e}^{2\alpha c}\!-\!1\big]\Big\}, 
 \end{alignat*}
with
\begin{equation*}
\widehat\gamma_j(u):=\exp \Big(2\beta \gamma_j\big(\log
(u+1)/2\alpha\big)\Big)-1, \qquad j\in \big\{\{1,2\},\{1\},\{2\},\{0\}\big\},
\end{equation*}
and obtain the following result.

\begin{Thm}
   \label{mainOU2} \ 
If \ $\alpha\not=0$ \ and \ $\beta\not=0$, \ functions \ $h_1, \ldots,
h_p$ \ are twice continuously differentiable   
inside \  $G$ \ and the partial derivatives  \ $\partial_1h_i, \ 
\partial_2h_i$ \ and \ $\partial_1\partial_2h_i, \ i=1, \ldots, p$, \ 
 can be continuously extended to  $G$ then the
probability measures 
 \ ${\mathsf P}_Z$ \ and \ ${\mathsf P}_X$, \ generated
 on \ $C(G)$ \ by the sheets \ $Z$ \ and 
 \ $X$, \ respectively, are equivalent and the
 Radon-Nikodym derivative of  
 \ ${\mathsf P}_Z$ \ with respect to \ ${\mathsf P}_X$ \ equals
 \begin{equation*}
  \frac{{\mathrm d}{\mathsf P}_Z}{{\mathrm d}{\mathsf
      P}_X}(Z)=\exp\left\{-\frac {\alpha\beta}{2\sigma
      ^2}\big({\mathbf m}^{\top}A{\mathbf 
      m}-2\zeta ^{\top} {\mathbf m}\big)\right\}, 
 \end{equation*}
where \ $A:=\big (A_{k,\ell} \big)_{k,\ell=1}^p$,  \ ${\mathbf
  m}:=(m_1,\ldots, m_p)^{\top}$ \ and \ 
  $\zeta :=\big(\zeta_1, \ldots , \zeta_p\big)^{\top}$ \ with
\begin{align}
  \label{eq:eq2.10}
A_{k,\ell}:\!&= \coth \big(\alpha a\big)\coth\big(\beta\gamma_2(a)\big)
h_k\big(a,\gamma_2(a)\big)h_{\ell}\big(a,\gamma_2(a)\big) 
+ h_k\big(c,\gamma_1(c)\big)h_{\ell}\big(c,\gamma_1(c)\big)
 \\ 
&+\coth\big(\beta\gamma_1(b_1)\big)h_k\big(b_1,\gamma_1(b_1)\big)
 h_{\ell}\big(b_1,\gamma_1(b_1)\big)
+ \coth \big(\alpha b_2\big)h_k\big(b_2,\gamma_2(b_2)\big)
h_{\ell}\big(b_2,\gamma_2(b_2)\big)\nonumber \\
&+\!\int\limits_a^{b_1}\!\!\coth\big(\beta\gamma_{1,2}(s)\big)\!\Big[\alpha
 h_k\big(s,\gamma_{1,2}(s)\big)h_{\ell}\big(s,\gamma_{1,2}(s)\big)
 \!+\!\alpha^{-1}\!\partial_1h_k\big(s,\gamma_{1,2}(s)\big) 
\partial_1h_{\ell}\big(s,\gamma_{1,2}(s)\big)\Big]{\mathrm d}s \nonumber \\ 
&+\int\limits_{b_1}^c\!\!\coth\big(\beta\gamma_1(s)\big)\Big[\alpha
 h_k\big(s,\gamma_1(s)\big)\!+\!\partial_1h_k\big(s,\gamma_1(s)\big)\Big]\Big[
 h_{\ell}\big(s,\gamma_1(s)\big)\!+\!\alpha^{-1}\partial_1
 h_{\ell}\big(s,\gamma_1(s)\big)\Big]{\mathrm d}s \nonumber \\ 
 &+\int\limits_a^{b_2}\Big[\alpha \coth\big(\alpha s\big)
 h_k\big(s,\gamma_2(s)\big)-\partial_1
 h_k\big(s,\gamma_2(s)\big)\Big] \nonumber\\
&\phantom{+\int\limits_a^{b_2}}\times
 \Big [\coth\big(\alpha s\big)h_{\ell}\big(s,\gamma_2(s)\big) 
 -\alpha^{-1}\partial_1h_{\ell}\big(s,\gamma_2(s)\big)\Big]{\mathrm
   d}s \nonumber 
\end{align}
\begin{align}
\phantom{A_{k,\ell}:\!}
&+\int\limits_{b_2}^c\Big[\alpha h_k\big(s,\gamma_0(s)\big)
h_{\ell}\big(s,\gamma_0(s)\big)+
\alpha^{-1}\partial_1h_k\big(s,\gamma_0(s)\big)
\partial_1h_{\ell}\big(s,\gamma_0(s)\big)\Big]{\mathrm d}s \nonumber
\\
&+\int\limits_{\gamma_{1,2}(b_1)}^{\gamma_{1,2}(a)}\coth\big(\alpha
\gamma_{1,2}^{-1}(t)\big)\Big[\beta\coth \big(\beta t\big) 
 h_k\big(\gamma_{1,2}^{-1}(t),t\big)-\partial_2
 h_k\big(\gamma_{1,2}^{-1}(t),t\big)\Big]  \nonumber \\ 
 &\phantom{===========.}
\times\Big[
\coth \big(\beta t\big)h_{\ell}\big(\gamma_{1,2}^{-1}(t),t\big)
-\beta^{-1}\partial_2 h_{\ell}\big(\gamma_{1,2}^{-1}(t),t\big)\Big]
{\mathrm d}t \nonumber \\ 
 &+\int\limits_{\gamma_1(b_1)}^{\gamma_1(c)}\Big[\beta
 h_k\big(\gamma_1^{-1}(t),t\big)h_{\ell}\big(\gamma_1^{-1}(t),t\big)
 +\beta^{-1}\partial_2h_k\big(\gamma_1^{-1}(t),t\big)
  \partial_2h_{\ell}\big(\gamma_1^{-1}(t),t\big)\!\Big]{\mathrm d}t
  \nonumber \\
 &+\!\!\int\limits_{\gamma_2(a)}^{\gamma_2(b_2)}\!\!\!\!\coth\big(\alpha
\gamma_2^{-1}(t)\big)\!\Big[\beta
 h_k\big(\gamma_2^{-1}(t),t\big)h_{\ell}\big(\gamma_2^{-1}(t),t\big)
 \!+\!\beta^{-1} \partial_2h_k\big(\gamma_2^{-1}(t),t\big) 
\partial_2h_{\ell}\big(\gamma_2^{-1}(t),t\big)\Big]{\mathrm d}t \nonumber \\
&+\int\limits_{\gamma_0(c)}^{\gamma_0(b_2)} \Big[\beta 
 h_k\big(\gamma_0^{-1}(t),t\big)+\partial_2
 h_k\big(\gamma_0^{-1}(t),t\big) \Big]\Big[
 h_{\ell}\big(\gamma_0^{-1}(t),t\big)+\beta^{-1}
\partial_2h_{\ell}\big(\gamma_0^{-1}(t),t\big)\Big]{\mathrm d}t \nonumber \\
&+\iint\limits_G \big [\alpha\beta
h_k(s,t)h_{\ell}(s,t)+\alpha^{-1}\beta \partial_1h_k(s,t)\partial_1h_{\ell}(s,t)
+\alpha\beta^{-1} \partial_2h_k(s,t)\partial_2h_{\ell}(s,t) \nonumber \\
&\phantom{======}
+\alpha^{-1}\beta^{-1} \partial_1\partial_2h_k(s,t)
\partial_1\partial_2h_{\ell}(s,t) \big]{\mathrm d}s\, {\mathrm d}t  \nonumber
\end{align}
and
\begin{align}
\label{eq:eq2.11}
\zeta_k:=&\,\Big[1+\coth (\alpha b_1)\Big]\Big[1+ \coth
\big(\beta\gamma_1(b_1)\big)\Big]
h_k\big(b_1,\gamma_1(b_1)\big)Z\big(b_1,\gamma_1(b_1)\big) \\ 
&+\int\limits_{b_1}^c\Big[1+ \coth
\big(\beta\gamma_1(s)\big)\Big]\Big[\alpha h_k\big(s,\gamma_1(s)\big)
 +\partial_1h_k\big(s,\gamma_1(s)\big)\Big] \nonumber \\
&\phantom{====} \times\Big[Z\big(s,\gamma_1(s)\big){\mathrm
  d}s+\alpha^{-1}Z\big({\mathrm 
  d}s,\gamma_1(s)\big)\Big]\nonumber\\
&+\int\limits_a^{b_1}\Big[1+ \coth
\big(\beta\gamma_{1,2}(s)\big)\Big]\Big[\alpha \coth (\alpha
s)h_k\big(s,\gamma_{1,2}(s)\big) 
 -\partial_1h_k\big(s,\gamma_{1,2}(s)\big)\Big] \nonumber \\
&\phantom{====}\times 
\Big[\coth (\alpha s)Z\big(s,\gamma_{1,2}(s)\big){\mathrm d}s
-\alpha^{-1}Z\big({\mathrm 
  d}s,\gamma_{1,2}(s)\big)\Big]\nonumber\\
&+\int\limits_{\gamma_2(a)}^{\gamma_2(b_2)}\Big[1+ \coth
\big(\alpha\gamma_2^{-1}(t)\big)\Big]\Big[\beta
 h_k\big(\gamma_2^{-1}(t),t\big)+ \partial_2h_k\big(\gamma_2^{-1}(t),t\big)\Big]
 \nonumber \\ 
&\phantom{====}\times 
\Big[Z\big(\gamma_2^{-1}(t),t\big){\mathrm d}t+
\beta^{-1}Z\big(\gamma_2^{-1}(t),{\mathrm d}t\big)\Big]\nonumber
\end{align}
\begin{align}
\phantom{\zeta_k:=}
&+\int\limits_{\gamma_{1,2}(b_1)}^{\gamma_{1,2}(a)}\Big[1+ \coth
\big(\alpha\gamma_{1,2}^{-1}(t)\big)\Big]\Big[\beta
 h_k\big(\gamma_{1,2}^{-1}(t),t\big)+ \partial_2h_k
 \big(\gamma_{1,2}^{-1}(t),t\big)\Big]  \nonumber \\
&\phantom{====}\times 
\Big[Z\big(\gamma_{1,2}^{-1}(t),t\big){\mathrm d}t+
\beta^{-1}Z\big(\gamma_{1,2}^{-1}(t),{\mathrm d}t\big)\Big] \nonumber\\ 
&+\iint\limits_G \big [\alpha\beta
h_k(s,t)+\beta \partial_1h_k(s,t)+\alpha\partial_2h_k(s,t)+
\partial_1\partial_2h_k(s,t)\big] \nonumber\\
&\phantom{====}\times \big [Z(s,t){\mathrm d}s\, {\mathrm d}t
+\alpha^{-1} Z({\mathrm d}s,t){\mathrm d}t  
+\beta^{-1} Z(s,{\mathrm d}t) {\mathrm d}s 
+\alpha^{-1}\beta^{-1} Z({\mathrm d}s,{\mathrm d}t)\big].  \nonumber
\end{align}
The maximum likelihood estimator of the parameter vector \ ${\mathbf
  m}$ \ based on the 
 observations \ $\{Z(s,t):(s,t)\in G\}$ \ has the form \ ${\widehat
   {\mathbf m}}=A^{-1}\zeta$ \ and 
 has a $p$-dimensional normal distribution with mean \ ${\mathbf m}$ \ and
 covariance matrix \ $A^{-1}$.
\end{Thm}

\begin{Rem}
  \label{compare} \
Observe that Theorems \ref{mainW},
\ref{mainOU1} and  \ref{mainOU2} are generalizations
 of Theorems 4, 5 and 6 of \citet{bpz2}, respectively, where 
 \ $G=[S_1,S_2]\times[T_1,T_2]$, \ with
 \ $[S_1,S_2],\ [T_1,T_2]\subset(0,\infty)$. \  Hence, from these
 theorems one can also derive the results of \citet{aratonm1,aratonm2}. 
\end{Rem}

\begin{figure}[tb]
 \begin{center}
\leavevmode
\epsfig{file=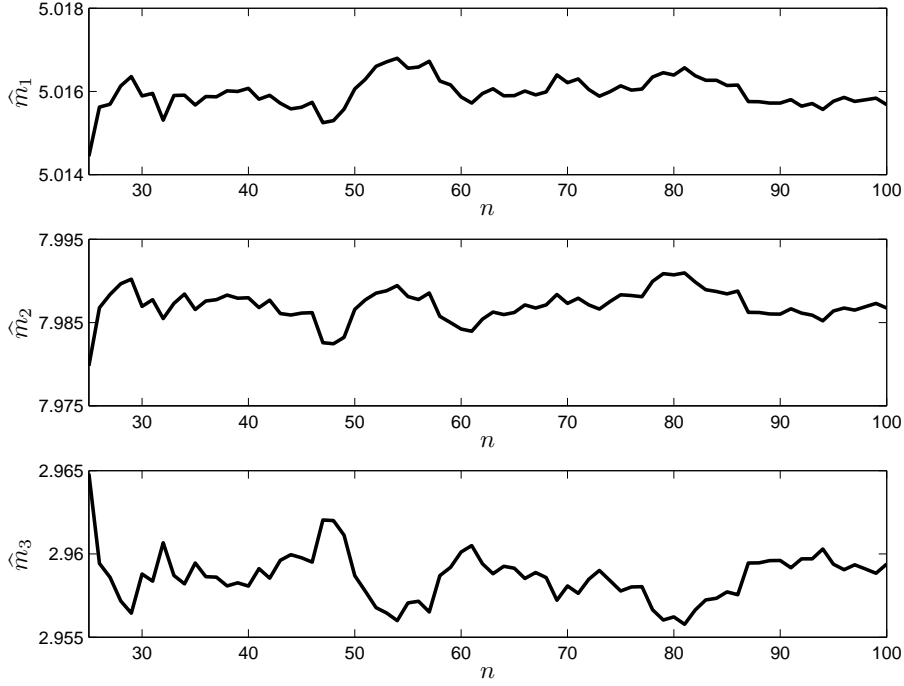,height=9.2cm}  
\caption{\label{fig:fig2} Means of the estimates of the components of
  \ ${\mathbf m}$ \  in Example \ref{ex1} for \ $25\leq n\leq 100$.}  
  \end{center}  
\end{figure}

\section{Simulation results}
  \label{sec:sec3}

\begin{figure}[t!]
 \begin{center}
\leavevmode
\epsfig{file=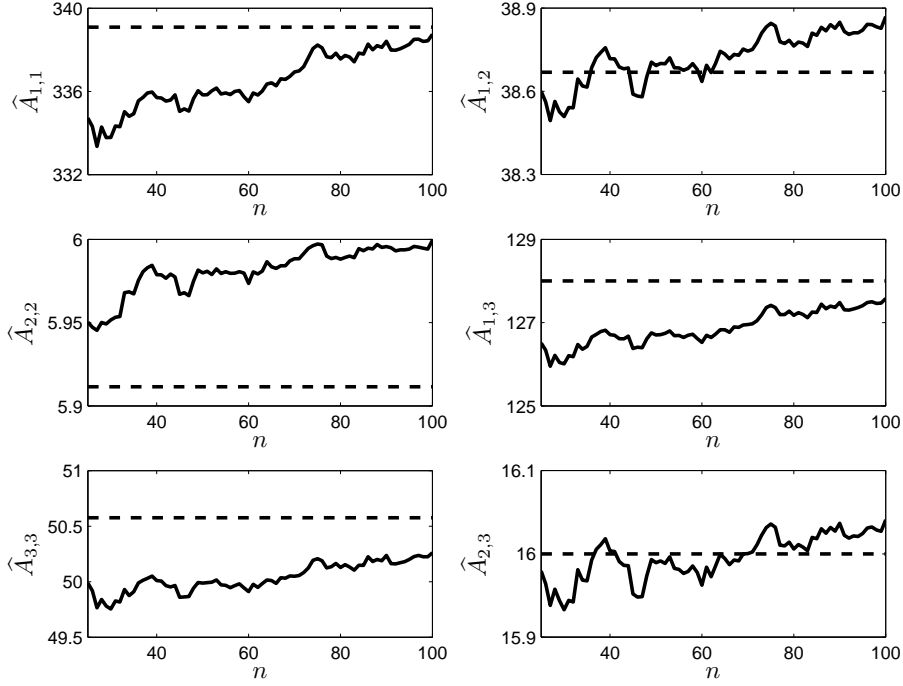,height=9.2cm}  
\caption{\label{fig:fig3} Estimated covariances of
  \ $\zeta$ \ in Example \ref{ex1} for \ $25\leq n\leq 100$.}  
  \end{center}  
\end{figure}
To illustrate our theoretical results we performed computer
simulations using Matlab 2010a.
In order to simulate the Gaussian random fields considered above their
Karhunen-Lo\`eve expansions are applied. For the Wiener sheet \ $W(s,t), 
\ 0\leq s \leq S, \ 0\leq t \leq T$ \ we have  
\begin{equation}
    \label{eq:eq3.1} 
 W(s,t)\approx  \sum\limits_{j,k=1}^{n} \omega_{j,k}
 \frac{8\sqrt{ST}}{\pi^2(2k-1)(2j-1)} \sin
 \left(\frac{\pi(2j-1)t}{2T}\right) 
 \sin \left(\frac{\pi(2k-1)s}{2S}\right),
 \end{equation}
where \ $\{\omega_{j,k}: 1\leq j,k\leq n\}$ \ are independent standard
normal random variables \citep{dpy}. The expansions for stationary and
zero start Ornstein-Uhlenbeck sheets \ $\widetilde X(s,t)$ \ and \
$X(s,t)$ \ with  \ $0\leq s \leq S, \ 0\leq t \leq T$ \ can directly
be derived from  
\eqref{eq:eq3.1} using representations \eqref{eq:eq2.4} and
\eqref{eq:eq2.8}, respectively (see e.g. \citet{jb}), yielding
\begin{align*}
 \widetilde X(s,t)\approx & \sum\limits_{j,k=1}^n \omega_{j,k}
 \frac{4 \sigma {\mathrm e}^{\alpha (S-s)+\beta (T-t)}}{\pi^2\sqrt
   {\alpha\beta }(2k-1)(2j-1)} \sin
 \left(\frac{\pi(2j-1) {\mathrm e}^{2\beta(t-T)}}2\right) 
 \sin \left(\frac{\pi(2k-1) {\mathrm e}^{2\alpha (s-S)}}2\right), \\
X(s,t)\approx & \sum\limits_{j,k=1}^n \omega_{j,k}
 \frac{4 \sigma \sqrt{({\mathrm e}^{2\alpha S}-1)({\mathrm e}^{2\beta
       T}\!-\!1)}}{\pi^2{\mathrm e}^{\alpha s+\beta t} \sqrt 
   {\alpha\beta } (2k-1)(2j-1)}
\sin
 \left(\frac{\pi(2j-1) ({\mathrm e}^{2\beta t}-1)}{2({\mathrm
       e}^{2\beta T}-1)}\right)  \\
 &\phantom{=====================}\times
 \sin \left(\frac{\pi(2k-1) ({\mathrm e}^{2\alpha s}-1)}{2({\mathrm
       e}^{2\alpha S}-1)}\right). 
 \end{align*}

In each of the following examples \ $1000$ \ independent samples of
the driving Gaussian sheet were simulated with $n$ varying between
$25$ and $100$ and the means of the estimates of the parameter vector
\ ${\mathbf m}$ \ and the empirical covariance matrices of the vectors
\ $\zeta$ \ defined by \eqref{eq:eq2.3}, \eqref{eq:eq2.7} and
\eqref{eq:eq2.11}, respectively, were calculated. 

\begin{Ex}
   \label{ex1} \
Consider the model 
\begin{equation*}
Z(s,t)=m_1(s^2+t^2)+m_2(s+t)+m_3(s\cdot t)+W(s,t),  \qquad (s,t)\in G,
\end{equation*}
where \ $W(s,t), \ (s,t)\in [0,8]^2$, \ is a standard Wiener sheet and
\ $G$ \ is a circle 
with center at \ $(6,6)$ \ and radius \ $r=2$. \ In this case the
entries of the matrix \ $A$ \ defined by \eqref{eq:eq2.2} and 
the approximations of the components of \
$\zeta=\big(\zeta_1,\zeta_2,\zeta_3)^{\top}$  \ defined by
\eqref{eq:eq2.3} can be calculated using numerical integration, where
Matlab function {\tt quad} is applied (recursive adaptive Simpson quadrature). 

\begin{figure}[tb]
 \begin{center}
\leavevmode
\epsfig{file=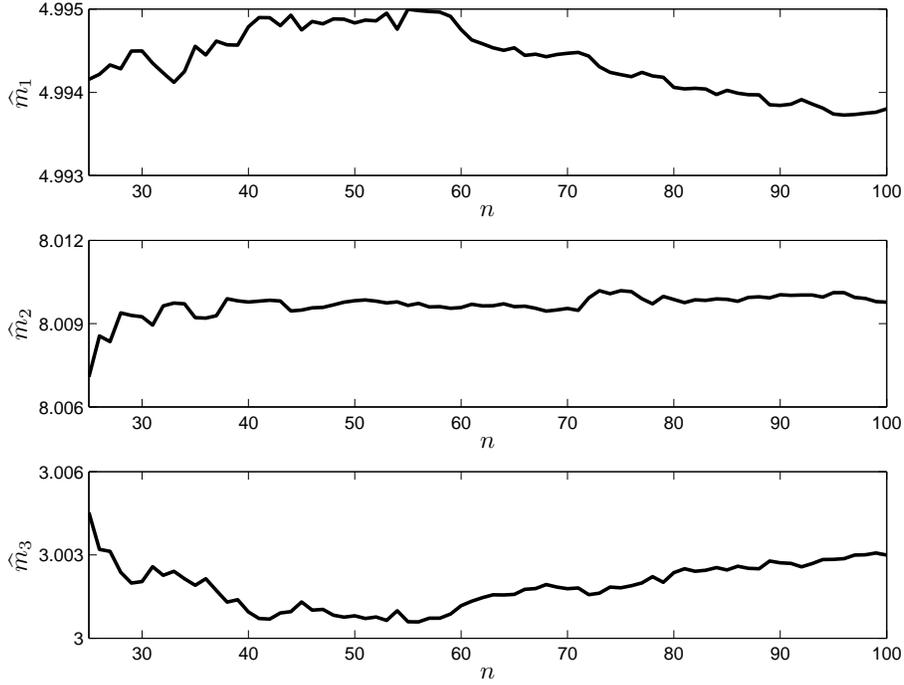,height=9.2cm}  
\caption{\label{fig:fig4} Means of the estimates of the components of
  \ ${\mathbf m}$ \  in Example \ref{ex2} for \ $25\leq n\leq 100$.}  
  \end{center}  
\end{figure}
The theoretical parameter values are \ $m_1=5, \ m_2=8$ \ and \
$m_3=3$, \ while the theoretical covariance matrix of \ $\zeta$ \ equals
 \begin{align*}
A&=  
  \begin{pmatrix}
  339.0895&38.6688&128.0000\\
   38.6688&5.9115&16.0000\\
   128.0000&16.0000&50.5752
  \end{pmatrix}.
 \end{align*} 
On Figure \ref{fig:fig2} the means of the estimates of the three
parameters, while on Figure \ref{fig:fig3} the estimated covariances of
\ $\zeta$  \ are plotted versus
the rate \ $n$ \ of the approximation \eqref{eq:eq3.1}. In case of \
$n=100$ \
we have \ $(5.0157,\ 7.9868,\ 2.9594)$ \ for the mean and
\begin{align*}
\widehat{A}&=  
  \begin{pmatrix}
  338.7473&38.8697&127.5784\\
   38.8697&5.9995&16.0409\\
  127.5784&16.0409&50.2610
  \end{pmatrix}
 \end{align*} 
for the covariance matrix. 
\end{Ex}

\begin{figure}[tb]
 \begin{center}
\leavevmode
\epsfig{file=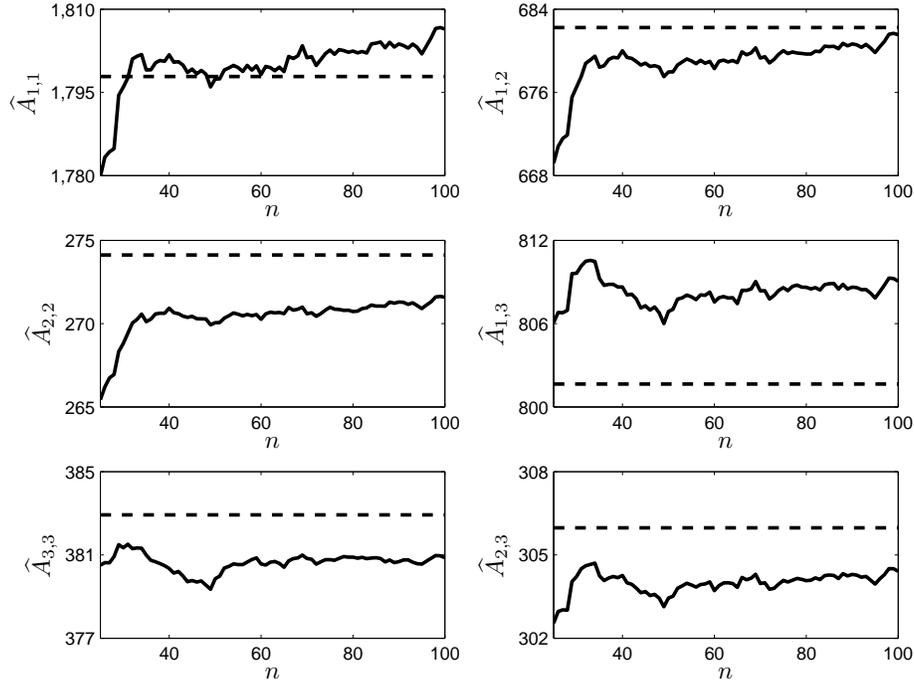,height=9.2cm}  
\caption{\label{fig:fig5} Estimated covariances of
  \ $\zeta$ \ in Example \ref{ex2} for \  $25\leq n\leq 100$.}  
  \end{center}  
\end{figure}

\begin{Ex}
   \label{ex2} \
Consider the model 
\begin{equation*}
Z(s,t)=m_1(s^2+t^2)+m_2(s+t)+m_3(s\cdot t)+\widetilde X(s,t),  \qquad (s,t)\in G,
\end{equation*}
where \ $\widetilde X(s,t), \ (s,t)\in [0,3]^2$, \ is a stationary
Ornstein-Uhlenbeck sheet with parameters $\alpha=1, \ \beta=1, \ \sigma=1$ and
\ $G$ \ is a circle 
with center at \ $(2,2)$ \ and radius \ $r=1$. \ Similarly to Example
\ref{ex1} the
entries of the matrix \ $A$ \ defined by \eqref{eq:eq2.6} and 
the approximations of the components of \
$\zeta=\big(\zeta_1,\zeta_2,\zeta_3)^{\top}$  \ defined by
\eqref{eq:eq2.7} can be calculated using numerical integration, where
Matlab functions {\tt quad} and  {\tt quad2d} \citep{shamp} are applied.  

The theoretical parameter values are the same as before, that is  \
$m_1=5, \ m_2=8$ \ and \ 
$m_3=3$, \ while the theoretical covariance matrix of \ $\zeta$ \
equals
 \begin{align*}
A&=  
  \begin{pmatrix}
  1797.8554&682.2301&801.6460\\
   682.2301&274.1195&305.9734\\
   801.6460&305.9734&382.9247
  \end{pmatrix}.
 \end{align*} 
On Figure \ref{fig:fig4} the means of the estimates of the three
parameters, while on Figure \ref{fig:fig5} the estimated covariances of
\ $\zeta$  \ are plotted versus
the rate \ $n$ \ of the Karhunen-Lo\`eve approximation. In case of \
$n=100$ \
we have \ $(4.9938,\ 8.0098,\ 3.0030)$ \ for the mean and
\begin{align*}
\widehat{A}&=  
  \begin{pmatrix}
  1806.4240&681.5420&809.0529\\
   681.5420&271.5779&304.4110\\
  809.0529&304.4110&380.8634
  \end{pmatrix}
 \end{align*} 
for the covariance matrix. 
\end{Ex}

\begin{figure}[tb]
 \begin{center}
\leavevmode
\epsfig{file=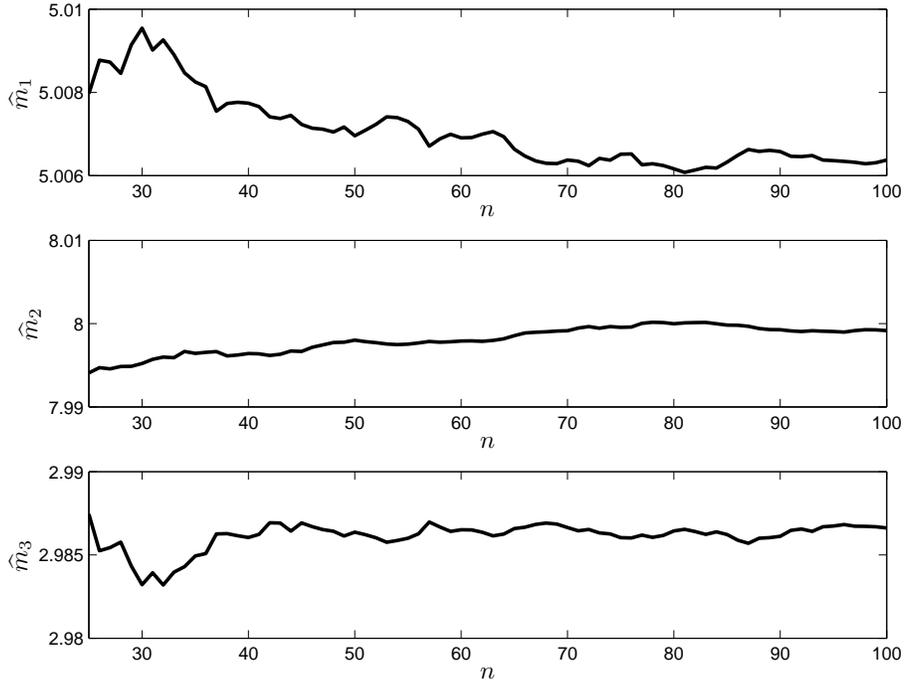,height=9.2cm}  
\caption{\label{fig:fig6} Means of the estimates of the components of
  \ ${\mathbf m}$ \ in Example \ref{ex3} for \ $25\leq n\leq 100$.}  
  \end{center}  
\end{figure}
\begin{figure}[tb]
 \begin{center}
\leavevmode
\epsfig{file=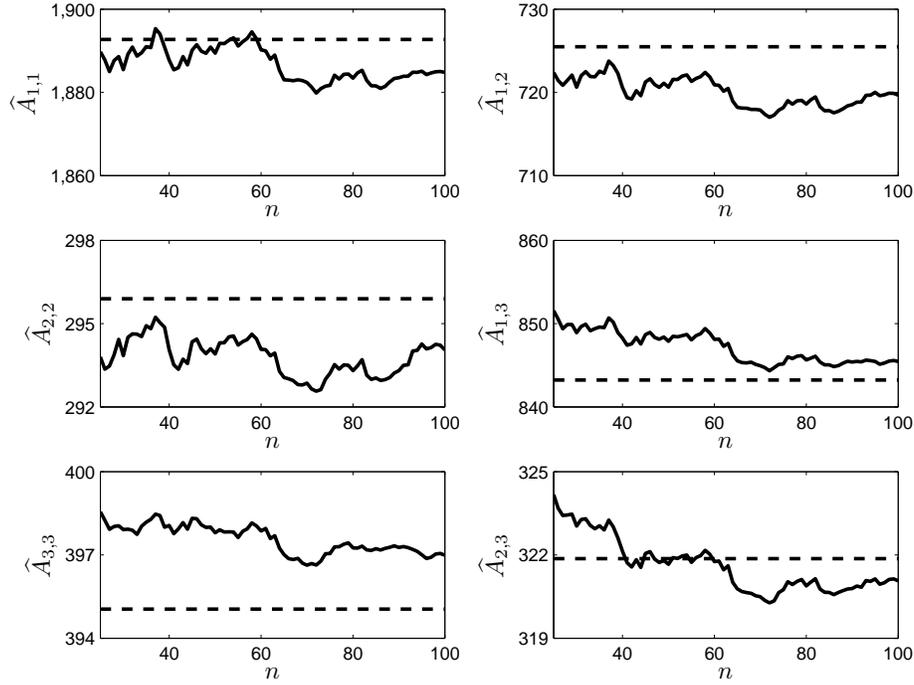,height=9.2cm}  
\caption{\label{fig:fig7} Estimated covariances of
  \ $\zeta$ \ in Example \ref{ex3} for \ $25\leq n\leq 100$.}  
  \end{center}  
\end{figure}
\begin{Ex}
   \label{ex3} \ Consider the same regression
\begin{equation*}
Z(s,t)=m_1(s^2+t^2)+m_2(s+t)+m_3(s\cdot t)+X(s,t),  \qquad (s,t)\in G,
\end{equation*}
as in Examples \ref{ex1} and \ref{ex2}, but now the driving process \
$X(s,t), \ (s,t)\in 
[0,3]^2$, \ is a zero start Ornstein-Uhlenbeck sheet with parameters
$\alpha=1, \ \beta=1, \ \sigma=1$. \ Similarly to Example
\ref{ex2}, \ $G$ \ is a circle 
with center at \ $(2,2)$ \ and radius \ $r=1$. \ Again, the
entries of the matrix \ $A$ \ defined by \eqref{eq:eq2.10} and 
the approximations of the components of \
$\zeta=\big(\zeta_1,\zeta_2,\zeta_3)^{\top}$  \ defined by
\eqref{eq:eq2.11} can only be calculated with the help of numerical integration.

The theoretical parameter values are the same as before, that is  \
$m_1=5, \ m_2=8$ \ and \ 
$m_3=3$, \ while the theoretical covariance matrix of \ $\zeta$ \
equals
 \begin{align*}
A&=  
  \begin{pmatrix}
  1892.7035&725.4822&843.2301\\
   725.4822&295.8952&321.8680\\
   843.2301&321.8680&395.0477
  \end{pmatrix}.
 \end{align*} 
Similarly to the previous examples, on Figure \ref{fig:fig6} the means
of the estimates of the three 
parameters, while on Figure \ref{fig:fig7} the estimated covariances of
\ $\zeta$  \ are plotted versus
the rate \ $n$ \ of the Karhunen-Lo\`eve approximation. In case of \
$n=100$ \
we have \ $(5.0067,\ 7.9992,\ 2.9866)$ \ for the mean and
\begin{align*}
\widehat{A}&=  
  \begin{pmatrix}
  1884.8293&719.6687&845.4606\\
   719.6687&294.0491&321.0663\\
  840.0300&321.0663&396.9913
  \end{pmatrix} 
 \end{align*} 
for the covariance matrix. 
\end{Ex}

\appendix
\section{Appendix}
   \label{sec:sec4}

\subsection{Proof of Theorem \ref{mainOU1}}
   \label{subsec4.1}

To prove Theorem \ref{mainOU1} one has to use representation
\eqref{eq:eq2.4} and apply Theorem \ref{mainW} 
for the random field  
\begin{equation*}
Y(u,v):=\frac {2\sqrt{\alpha\beta uv}}\sigma Z\Big(\frac{\log
  u}{2\alpha}, \frac{\log v}{2\beta}\Big)=m_1g_1(u,v)+\ldots
+m_pg_p(u,v)+W(u,v)
\end{equation*}
observed on \ $\widetilde G$, \ where functions \ $g_k$ \ are defined
by \eqref{eq:eq2.5}. As 
\begin{align}
\partial_1g_k(u,v)&=\frac {\sqrt{\beta v}}{\sigma \sqrt{\alpha
    u}}(\alpha+\partial_1)h_k\Big(\frac{\log u}{2\alpha}, \frac{\log
  v}{2\beta}\Big),  \nonumber \\
\partial_2g_k(u,v)&=\frac {\sqrt{\alpha u}}{\sigma \sqrt{\beta
    v}}(\beta+\partial_2)h_k\Big(\frac{\log u}{2\alpha}, \frac{\log
  v}{2\beta}\Big),  \label{eq:eq4.1}\\
\partial_1\partial_2g_k(u,v)&=\frac 1{2\sigma \sqrt{\alpha\beta
    uv}}(\alpha+\partial_1)(\beta+\partial_2)h_k\Big(\frac{\log
  u}{2\alpha}, \frac{\log  
  v}{2\beta}\Big),\nonumber
\end{align}
short calculations show
\begin{align*}
A_{k,\ell}=& \,4h_k\big(b_1,\gamma_{1,2}(b_1)\big)\,
    h_{\ell}\big(b_1,\gamma_{1,2}(b_1)\big)\\
&+2\int\limits_a^{b_1}\Big[(\alpha -\partial_1)
h_k\big(s,\gamma_{1,2}(s)\big)\Big]\Big[(1 -\alpha^{-1}\partial_1)
h_{\ell}\big(s,\gamma_{1,2}(s)\big)\Big] \,{\mathrm d}s\\
&+2\int\limits_{b_1}^c\Big[(\alpha +\partial_1)
h_k\big(s,\gamma_1(s)\big)\Big]\Big[(1 +\alpha^{-1}\partial_1)
h_{\ell}\big(s,\gamma_1(s)\big)\Big] \,{\mathrm d}s\\
&+2\int\limits_{\gamma_2(a)}^{\gamma_2(b_2)} \Big [(\beta+\partial_2)
  h_k\big(\gamma_2^{-1}(t),t\big)\Big]\Big [(1+\beta^{-1}\partial_2)
  h_{\ell}\big(\gamma_2^{-1}(t),t\big)\Big]\,{\mathrm d} t \\
&+2\int\limits_{\gamma_{1,2}(b_1)}^{\gamma_{1,2}(a)} \Big [(\beta+\partial_2)
  h_k\big(\gamma_{1,2}^{-1}(t),t\big)\Big]\Big[(1+\beta^{-1}\partial_2)
  h_{\ell}\big(\gamma_{1,2}^{-1}(t),t\big)\Big]\,{\mathrm d} t \\
&+\iint\limits_G\Big[(\alpha+\partial_1)(\beta+\partial_2)h_k(s,t)
\Big]\Big[(1+\alpha^{-1}\partial_1)(1+\beta^{-1}\partial_2) 
h_{\ell}(s,t)\Big]\,{\mathrm d}s\,{\mathrm d} t\\
=&A_{k,\ell}^{(1)}+\alpha A_{k,\ell}^{(2)}+\beta
A_{k,\ell}^{(3)}+\alpha^{-1}
A_{k,\ell}^{(4)}+\beta^{-1}A_{k,\ell}^{(5)}+A_{k,\ell}^{(6)},
\end{align*}
where
\begin{align*}
A_{k,\ell}^{(1)}:=&\,4h_k\big(b_1,\gamma_{1,2}(b_1)
h_{\ell}\big(b_1,\gamma_{1,2}(b_1)\big)+ 
2\int\limits_{b_1}^c\partial_1\Big[h_k\big(s,\gamma_1(s)\big)
h_{\ell}\big(s,\gamma_1(s)\big) \Big]\,{\mathrm d}s \\
&- 2\int\limits_a^{b_1}\partial_1\Big[h_k\big(s,\gamma_{1,2}(s)\big)
h_{\ell}\big(s,\gamma_{1,2}(s)\big)\Big]\,{\mathrm d}s
+2\int\limits_{\gamma_2(a)}^{\gamma_2(b_2)} \partial_2\Big[
  h_k\big(\gamma_2^{-1}(t),t\big)h_{\ell}
  \big(\gamma_2^{-1}(t),t\big)\Big]\,{\mathrm d} t\\  
&+2\int\limits_{\gamma_{1,2}(b_1)}^{\gamma_{1,2}(a)} \partial_2\Big [
  h_k\big(\gamma_{1,2}^{-1}(t),t\big)h_{\ell}
  \big(\gamma_{1,2}^{-1}(t),t\big)\Big]\,{\mathrm d} t
+\iint\limits_G \partial_1\partial_2
\big[h_k(s,t)h_{\ell}(s,t)\big]\,{\mathrm d}s\,{\mathrm d} t, \\
A_{k,\ell}^{(2)}:=&\,2\int\limits_{b_1}^ch_k\big(s,\gamma_1(s)\big)
h_{\ell}\big(s,\gamma_1(s)\big) \,{\mathrm d}s 
+2\int\limits_a^{b_1}h_k\big(s,\gamma_{1,2}(s)\big)
h_{\ell}\big(s,\gamma_{1,2}(s)\big)\,{\mathrm d}s \\
&+\iint\limits_G \partial_2
\big[h_k(s,t)h_{\ell}(s,t)\big]\,{\mathrm d}s\,{\mathrm d} t, \\
A_{k,\ell}^{(3)}:=&2\int\limits_{\gamma_2(a)}^{\gamma_2(b_2)}
h_k\big(\gamma_2^{-1}(t),t\big)h_{\ell}\big(\gamma_2^{-1}(t),t\big)\,{\mathrm
  d}t +2\int\limits_{\gamma_{1,2}(b_1)}^{\gamma_{1,2}(a)}
  h_k\big(\gamma_{1,2}^{-1}(t),t\big)
  h_{\ell}\big(\gamma_{1,2}^{-1}(t),t\big)\,{\mathrm d} t \\
&+\iint\limits_G \partial_1
\big[h_k(s,t)h_{\ell}(s,t)\big]\,{\mathrm d}s\,{\mathrm d} t, \\
A_{k,\ell}^{(4)}:=&\,2\int\limits_{b_1}^c\partial_1h_k\big(s,\gamma_1(s)\big)\,
\partial_1h_{\ell}\big(s,\gamma_1(s)\big) \,{\mathrm d}s 
+2\int\limits_a^{b_1}\partial_1h_k\big(s,\gamma_{1,2}(s)\big)\,
\partial_1h_{\ell}\big(s,\gamma_{1,2}(s)\big)\,{\mathrm d}s \\
&+\iint\limits_G \partial_2
\big[\partial_1h_k(s,t)\,\partial_1h_{\ell}(s,t)\big] \,{\mathrm
  d}s\,{\mathrm d} t, \\
A_{k,\ell}^{(5)}:=&2\int\limits_{\gamma_2(a)}^{\gamma_2(b_1)}
\partial_2h_k\big(\gamma_2^{-1}(t),t\big)\, 
\partial_2h_{\ell}\big(\gamma_2^{-1}(t),t\big)\,{\mathrm 
  d}t +2\int\limits_{\gamma_{1,2}(b_1)}^{\gamma_{1,2}(a)}
  \partial_2h_k\big(\gamma_{1,2}^{-1}(t),t\big)\,
  \partial_2h_{\ell}\big(\gamma_{1,2}^{-1}(t),t\big)\,{\mathrm d} t \\
&+\iint\limits_G \partial_1
\big[\partial_2h_k(s,t)\,\partial_2h_{\ell}(s,t)\big]\,{\mathrm
  d}s\,{\mathrm d} t, \\ 
A_{k,\ell}^{(6)}:=&\iint\limits_G \big [\alpha\beta
h_k(s,t)h_{\ell}(s,t)+\alpha^{-1}\beta \partial_1h_k(s,t)\partial_1h_{\ell}(s,t)
+\alpha\beta^{-1} \partial_2h_k(s,t)\partial_2h_{\ell}(s,t) \\
&\phantom{+\iint\limits_G=}
+\alpha^{-1}\beta^{-1} \partial_1\partial_2h_k(s,t)
\partial_1\partial_2h_{\ell}(s,t) \big]\,{\mathrm d}s\, {\mathrm d}t .
\end{align*}
Obviously
\begin{align}
\iint\limits_G \partial_2
\big[h_k(s,t)&\,h_{\ell}(s,t)\big]{\mathrm d}s\,{\mathrm d} t= \!
\int\limits_a^{b_2}\! h_k\big(s,\gamma_2(s)\big)
h_{\ell}\big(s,\gamma_2(s)\big)\,{\mathrm
  d}s-\! \int\limits_a^{b_1}\! h_k\big(s,\gamma_{1,2}(s)\big) 
h_{\ell}\big(s,\gamma_{1,2}(s)\big)\,{\mathrm d}s \nonumber\\ 
&-\! \int\limits_{b_1}^c \! h_k\big(s,\gamma_1(s)\big)
h_{\ell}\big(s,\gamma_1(s)\big)\,{\mathrm d}s
+\! \int\limits_{b_2}^c\! h_k\big(s,\gamma_0(s)\big)
h_{\ell}\big(s,\gamma_0(s)\big)\,{\mathrm d}s, \label{eq:eq4.2} \\
\iint\limits_G \partial_1
\big[h_k(s,t)&\,h_{\ell}(s,t)\big]{\mathrm d}s\,{\mathrm d} t=\!\!\!
\int\limits_{\gamma_1(b_1)}^{\gamma_1(c)}\!\!\!  h_k
\big(\gamma_1^{-1}(t),t\big)h_{\ell}\big(\gamma_1^{-1}(t),t\big)\,{\mathrm
  d}t -\!\!\!\!\!\!\int\limits_{\gamma_{1,2}(b_1)}^{\gamma_{1,2}(a)}\!\!\!\!\!
  h_k\big(\gamma_{1,2}^{-1}(t),t\big)
  h_{\ell}\big(\gamma_{1,2}^{-1}(t),t\big)\,{\mathrm d} t \nonumber\\
&-\!\!\int\limits_{\gamma_2(a)}^{\gamma_2(b_2)}\!\!
h_k\big(\gamma_2^{-1}(t),t\big)h_{\ell}\big(\gamma_2^{-1}(t),t\big)\,{\mathrm
  d}t +\!\!\!\int\limits_{\gamma_0(c)}^{\gamma_0(b_2)}\!\!
  h_k\big(\gamma_0^{-1}(t),t\big)
  h_{\ell}\big(\gamma_0^{-1}(t),t\big)\,{\mathrm d} t, \label{eq:eq4.3}
\end{align}
and similar expressions can be derived for
\begin{equation*}
\iint\limits_G \partial_2
\big[\partial_1h_k(s,t)\,\partial_1h_{\ell}(s,t)\big]\,{\mathrm
  d}s\,{\mathrm d} t  \qquad \text{and} \qquad
\iint\limits_G \partial_1
\big[\partial_2h_k(s,t)\,\partial_2h_{\ell}(s,t)\big] \,{\mathrm
  d}s\,{\mathrm d} t,  
\end{equation*}
respectively. Hence, using also that
\begin{align}
2\iint\limits_G &\partial_1\partial_2
\big[h_k(s,t)h_{\ell}(s,t)\big]\,{\mathrm d}s\,{\mathrm d} t \label{eq:eq4.4}\\ 
=& \int\limits_a^{b_2} \partial_1\Big[h_k\big(s,\gamma_2(s)\big)
h_{\ell}\big(s,\gamma_2(s)\big)\Big]\,{\mathrm d}s
-\int\limits_{\gamma_2(a)}^{\gamma_2(b_2)}
\partial_2\Big[h_k\big(\gamma_2^{-1}(t),t\big)
h_{\ell}\big(\gamma_2^{-1}(t),t\big)\Big] \,{\mathrm d}t \nonumber \\
&- \int\limits_{b_1}^c \partial_1\Big[h_k\big(s,\gamma_1(s)\big)
h_{\ell}\big(s,\gamma_1(s)\big)\Big]\,{\mathrm d}s
+\int\limits_{\gamma_1(b_1)}^{\gamma_1(c)}
\partial_2\Big[h_k\big(\gamma_1^{-1}(t),t\big)
h_{\ell}\big(\gamma_1^{-1}(t),t\big)\Big] \,{\mathrm d}t \nonumber \\
&- \int\limits_a^{b_1} \partial_1\Big[h_k\big(s,\gamma_{1,2}(s)\big)
h_{\ell}\big(s,\gamma_{1,2}(s)\big)\Big]\,{\mathrm d}s
+\int\limits_{\gamma_{1,2}(a)}^{\gamma_{1,2}(b_1)}
\partial_2\Big[h_k\big(\gamma_{1,2}^{-1}(t),t\big)
h_{\ell}\big(\gamma_{1,2}^{-1}(t),t\big)\Big] \,{\mathrm d}t \nonumber
\\
&+ \int\limits_{b_2}^c \partial_1\Big[h_k\big(s,\gamma_0(s)\big)
h_{\ell}\big(s,\gamma_0(s)\big)\Big]\,{\mathrm d}s
-\int\limits_{\gamma_0(b_2)}^{\gamma_0(c)}
\partial_2\Big[h_k\big(\gamma_0^{-1}(t),t\big)
h_{\ell}\big(\gamma_0^{-1}(t),t\big)\Big] \,{\mathrm d}t, \nonumber
\end{align}
after tedious but straightforward calculations we obtain
\eqref{eq:eq2.6}. 

Further, using again \eqref{eq:eq4.1} and Remark \ref{integral} we have
\begin{equation*}
\zeta_k=\zeta_k^{(1)}+\zeta_k^{(2)}+\zeta_k^{(3)}+\zeta_k^{(4)}+\zeta_k^{(5)}
+\zeta_k^{(6)},  
\end{equation*}
where
\begin{align*}
\zeta_k^{(1)}:=&\,4h_k\big(b_1,\gamma_{1,2}(b_1)\big)
Z\big(b_1,\gamma_{1,2}(b_1)\big)\\
\zeta_k^{(2)}:=&\, \qlim _{\varrho\to 0} \frac 4{\varrho}
\int\limits_{b_1}^c {\mathrm e}^{\alpha s}
(\alpha+\partial_1)h_k\big(s,\gamma_1(s)\big) \\
&\phantom{\, \qlim _{\varrho\to 0} \frac 1{\varrho}
\int\limits_{b_1}^c }
\times\bigg[\big({\mathrm e}^{2\alpha
  s}+\varrho\big)^{1/2}Z\Big(\frac{\log({\mathrm 
  e}^{2\alpha s}+\varrho)}{2\alpha},\gamma_1(s)\Big) -{\mathrm
  e}^{\alpha s}Z\big(s,\gamma_1(s)\big)\bigg]\,{\mathrm d}s, \\
\zeta_k^{(3)}:=&\, 4\int\limits_a^{b_1} \big(\alpha
-\partial_1\big)h_k\big(s,\gamma_{1,2}(s)\big)Z\big(s,\gamma_{1,2}(s)\big)
{\mathrm d}s \\
&-\qlim _{\varrho\to 0} \frac 4{\varrho}
\int\limits_a^{b_1} {\mathrm e}^{\alpha s}
(\alpha-\partial_1)h_k\big(s,\gamma_{1,2}(s)\big) \\
&\phantom{\, \qlim _{\varrho\to 0} \frac 1{\varrho}
\int\limits_a^{b_1} }
\times \bigg[\big({\mathrm e}^{2\alpha
  s}+\varrho\big)^{1/2}Z\Big(\frac{\log({\mathrm 
  e}^{2\alpha s}+\varrho)}{2\alpha},\gamma_{1,2}(s)\Big) -{\mathrm
  e}^{\alpha s}Z\big(s,\gamma_{1,2}(s)\big)\bigg]\,{\mathrm d}s, \\
\zeta_k^{(4)}:=&\, \qlim _{\delta\to 0} \frac 4{\delta}
\int\limits_{\gamma_2(a)}^{\gamma_2(b_2)} {\mathrm e}^{\beta t}
(\beta+\partial_2)h_k\big(\gamma_2^{-1}(t),t\big) \\
&\phantom{\, \qlim _{\varrho\to 0} \frac 1{\varrho}
\int\limits_{\gamma_2(a)}^{\gamma_2(b_2)} }
\times\bigg[\big({\mathrm e}^{2\beta
  t}+\delta\big)^{1/2}Z\Big(\gamma_2^{-1}(t),\frac{\log({\mathrm 
  e}^{2\beta t}+\delta)}{2\beta}\Big) -{\mathrm
  e}^{\beta t}Z\big(\gamma_2^{-1}(t),t\big)\bigg]\,{\mathrm d}t, \\
\zeta_k^{(5)}:=&\, \qlim _{\delta\to 0} \frac 4{\delta}
\int\limits_{\gamma_{1,2}(b_1)}^{\gamma_{1,2}(a)} {\mathrm e}^{\beta t}
(\beta+\partial_2)h_k\big(\gamma_{1,2}^{-1}(t),t\big) \\
&\phantom{\, \qlim _{\varrho\to 0} \frac 1{\varrho}
\int\limits_{\gamma_{1,2}(b_1)}^{\gamma_{1,2}(1)} }
\times\bigg[\big({\mathrm e}^{2\beta
  t}+\delta\big)^{1/2}Z\Big(\gamma_{1,2}^{-1}(t),\frac{\log({\mathrm 
  e}^{2\beta t}+\delta)}{2\beta}\Big) -{\mathrm
  e}^{\beta t}Z\big(\gamma_{1,2}^{-1}(t),t\big)\bigg]\,{\mathrm d}t, \\
\zeta_k^{(6)}:=&\,\qlim _{\varrho,\delta\to 0} \frac 4{\varrho \delta}
\iint\limits _G  {\mathrm e}^{\alpha s+\beta t}
(\alpha+\partial_1)(\beta+\partial_2)h_k\big(s,t\big)\bigg[{\mathrm e}^{\alpha
  s+\beta t}Z(s,t) \\
&\phantom{,\qlim _{\varrho,\delta\to 0} \frac 4{\varrho \delta}
\iint\limits _G }
-\!\big({\mathrm e}^{2\alpha
  s}\!+\!\varrho\big)^{1/2}{\mathrm e}^{\beta
  t}Z\Big(\frac{\log({\mathrm 
  e}^{2\alpha s}\!+\!\varrho)}{2\alpha},t\Big)\!-\!{\mathrm e}^{\alpha
  s}\big({\mathrm e}^{2\beta
  t}\!+\!\delta\big)^{1/2}Z\Big(s,\frac{\log({\mathrm 
  e}^{2\beta t}\!+\!\delta)}{2\beta}\Big)\\
&\phantom{,\qlim _{\varrho,\delta\to 0} \frac 4{\varrho \delta}
\iint\limits _G }
+\!\big({\mathrm e}^{2\alpha
  s}+\varrho\big)^{1/2}\big({\mathrm e}^{2\beta
  t}+\delta\big)^{1/2}Z\Big(\frac{\log({\mathrm 
  e}^{2\alpha s}+\varrho)}{2\alpha},\frac{\log({\mathrm 
  e}^{2\beta t}+\delta)}{2\beta}\Big)\bigg]{\mathrm d}s\,{\mathrm d}t,
\end{align*}
implying \eqref{eq:eq2.7}. \proofend

\subsection{Proof of Theorem \ref{mainOU2}}
   \label{subsec4.2}
The proof, which we give only in the case \ $\alpha >0, \ \beta>0$, \ is
similar to the proof of  Theorem \ref{mainOU1}. Here 
one has to use representation \eqref{eq:eq2.8} and apply Theorem \ref{mainW} 
for the random field  
\begin{align*}
Y(u,v):=&\,\frac {2\sqrt{\alpha\beta (u+1)(v+1)}}\sigma Z\Big(\frac{\log
  (u+1)}{2\alpha}, \frac{\log (v+1)}{2\beta}\Big) \\
=&\,m_1g_1(u,v)+\ldots
+m_pg_p(u,v)+W(u,v)
\end{align*}
observed on \ $\widehat G$, \ and functions \ $g_k$ \ should be defined
by \eqref{eq:eq2.9}. In this way
\begin{align}
\partial_1g_k(u,v)&=\frac {\sqrt{\beta (v+1)}}{\sigma \sqrt{\alpha
    (u+1)}}(\alpha+\partial_1)h_k\Big(\frac{\log (u+1)}{2\alpha}, \frac{\log
  (v+1)}{2\beta}\Big),  \nonumber \\
\partial_2g_k(u,v)&=\frac {\sqrt{\alpha (u+1)}}{\sigma \sqrt{\beta
    (v+1)}}(\beta+\partial_2)h_k\Big(\frac{\log (u+1)}{2\alpha}, \frac{\log
  (v+1)}{2\beta}\Big),  \label{eq:eq4.5}\\
\partial_1\partial_2g_k(u,v)&=\frac 1{2\sigma \sqrt{\alpha\beta
    (u+1)(v+1)}}(\alpha+\partial_1)(\beta+\partial_2)h_k\Big(\frac{\log
  (u+1)}{2\alpha}, \frac{\log  
  (v+1)}{2\beta}\Big),\nonumber
\end{align}
so
\begin{align*}
&A_{k,\ell}= \,4\frac {{\mathrm e}^{2\alpha b_1}}{{\mathrm
    e}^{2\alpha b_1}-1}\frac {{\mathrm e}^{2\beta \gamma_{1,2}(b_1)}}
{{\mathrm e}^{2\beta \gamma_{1,2}(b_1)}-1} 
h_k\big(b_1,\gamma_{1,2}(b_1)\big)\,
    h_{\ell}\big(b_1,\gamma_{1,2}(b_1)\big)\\
&\phantom{=}+2\!\int\limits_a^{b_1}\!\frac {{\mathrm e}^{2\beta \gamma_{1,2}(s)}}
{{\mathrm e}^{2\beta \gamma_{1,2}(s)}\!-\!1}
\Big[\big(\alpha\coth(\alpha s)\!-\!\partial_1\big) 
h_k\big(s,\gamma_{1,2}(s)\big)\Big]\Big[\big(\coth(\alpha
s)\!-\!\alpha^{-1}\partial_1\big) 
h_{\ell}\big(s,\gamma_{1,2}(s)\big)\Big] \,{\mathrm d}s\\
&\phantom{=}+2\int\limits_{b_1}^c\frac {{\mathrm e}^{2\beta \gamma_1(s)}}
{{\mathrm e}^{2\beta \gamma_1(s)}\!-\!1}\Big[(\alpha +\partial_1)
h_k\big(s,\gamma_1(s)\big)\Big]\Big[(1 +\alpha^{-1}\partial_1)
h_{\ell}\big(s,\gamma_1(s)\big)\Big] \,{\mathrm d}s\\
&\phantom{=}+2\int\limits_{\gamma_2(a)}^{\gamma_2(b_2)} \frac
{{\mathrm e}^{2\alpha \gamma_2^{-1}(t)}} 
{{\mathrm e}^{2\alpha \gamma_2^{-1}(t)}\!-\!1} \Big [(\beta+\partial_2)
  h_k\big(\gamma_2^{-1}(t),t\big)\Big]\Big[(1+\beta^{-1}\partial_2)
  h_{\ell}\big(\gamma_2^{-1}(t),t\big)\Big]\,{\mathrm d} t \\
&\phantom{=}+2\int\limits_{\gamma_{1,2}(b_1)}^{\gamma_{1,2}(a)}
\frac {{\mathrm e}^{2\alpha \gamma_{1,2}^{-1}(t)}}
{{\mathrm e}^{2\alpha \gamma_{1,2}^{-1}(t)}\!-\!1} \Big [(\beta+\partial_2)
  h_k\big(\gamma_{1,2}^{-1}(t),t\big)\Big]\Big[(1+\beta^{-1}\partial_2)
  h_{\ell}\big(\gamma_{1,2}^{-1}(t),t\big)\Big]\,{\mathrm d} t \\
&\phantom{=}+\iint\limits_G\Big[(\alpha+\partial_1)(\beta+\partial_2)h_k(s,t)
\Big]\Big[(1+\alpha^{-1}\partial_1)(1+\beta^{-1}\partial_2) 
h_{\ell}(s,t)\Big]\,{\mathrm d}s\,{\mathrm d} t.
\end{align*}
Collecting separately the terms containing \ $\gamma_{1,2}, \
\gamma_1, \ \gamma_2$ \ and \ $\gamma_0$, \ after long straightforward
calculations using again \eqref{eq:eq4.2}--\eqref{eq:eq4.4} we obtain
\begin{equation*}
A_{k,\ell}=A_{k,\ell}^{(1)}+A_{k,\ell}^{(2)}+A_{k,\ell}^{(3)}+A_{k,\ell}^{(4)}
+A_{k,\ell}^{(5)},
\end{equation*}
where
\begin{align*}
A_{k,\ell}^{(1)}:=& \, \frac 12  h_k\big(b_1,\gamma_{1,2}(b_1)\big) 
h_{\ell}\big(b_1,\gamma_{1,2}(b_1)\big)-\frac 12  h_k\big(a,\gamma_{1,2}(a)\big) 
h_{\ell}\big(a,\gamma_{1,2}(a)\big)\\
&+\coth\big(\alpha
a\big)\coth\big(\beta\gamma_{1,2}(a)\big)h_k\big(a,\gamma_{1,2}(a)\big)  
h_{\ell}\big(a,\gamma_{1,2}(a)\big) \\
&+\coth\big(\beta\gamma_{1,2}(b_1)\big)
h_k\big(b_1,\gamma_{1,2}(b_1)\big) 
h_{\ell}\big(b_1,\gamma_{1,2}(b_1)\big)+\coth\big(\alpha
a\big)h_k\big(a,\gamma_{1,2}(a)\big)  
h_{\ell}\big(a,\gamma_{1,2}(a)\big)\\
&+\int\limits_a^{b_1} \coth\big(\beta\gamma_{1,2}(s)\big) \Big[ \alpha
h_k\big(s,\gamma_{1,2}(s)\big) 
h_{\ell}\big(s,\gamma_{1,2}(s)\big)+\alpha^{-1} \partial_1
h_k\big(s,\gamma_{1,2}(s)\big) 
\partial _1h_{\ell}\big(s,\gamma_{1,2}(s)\big)\Big] {\mathrm d}s \\
&+\int\limits_{\gamma_{1,2}(b_1)}^{\gamma_{1,2}(a)}
\coth\big(\alpha\gamma_{1,2}^{-1}(t)\big) 
\Big[\big(\beta \coth (\beta t)-\partial_2\big)
h_k\big(\gamma_{1,2}^{-1}(t),t\big) \Big] \\
&\phantom{============}\times\Big[
\big(\coth (\beta t)-\beta^{-1}\partial _2\big)
h_{\ell}\big(\gamma_{1,2}^{-1}(t),t\big)\Big] {\mathrm d}t, \\ 
A_{k,\ell}^{(2)}:=& \, \frac 12  h_k\big(c,\gamma_1(c)\big) 
h_{\ell}\big(c,\gamma_1(c)\big)-\frac 12  h_k\big(b_1,\gamma_1(b_1)\big) 
h_{\ell}\big(b_1,\gamma_1(b_1)\big)\\
&+\int\limits_{b_1}^c \coth\big(\beta\gamma_1(s)\big) \Big[\big(\alpha
+\partial _1\big) h_k\big(s,\gamma_1(s)\big)\Big]\Big[\big(1+
\alpha^{-1}\partial_1 \big)h_{\ell}\big(s,\gamma_1(s)\big)\Big] {\mathrm d}s \\
&+\int\limits_{\gamma_1(b_1)}^{\gamma_1(c)}
\Big[\beta h_k\big(\gamma_1^{-1}(t),t\big)h_{\ell}\big(\gamma_1^{-1}(t),t\big)
+\beta^{-1}\partial_2h_k\big(\gamma_1^{-1}(t),t\big)\partial _2
h_{\ell}\big(\gamma_1^{-1}(t),t\big)\Big] {\mathrm d}t, \\
A_{k,\ell}^{(3)}:=& \, \frac 12  h_k\big(b_2,\gamma_2(b_2)\big) 
h_{\ell}\big(b_2,\gamma_2(b_2)\big)-\frac 12  h_k\big(a,\gamma_2(a)\big) 
h_{\ell}\big(a,\gamma_2(a)\big)\\
&+\coth\big(\alpha b_2\big)h_k\big(b_2,\gamma_2(b_2)\big)  
h_{\ell}\big(b_2,\gamma_2(b_2)\big)-\coth\big(\alpha a\big)
h_k\big(a,\gamma_2(a)\big) h_{\ell}\big(a,\gamma_2(a)\big)\\
&+\int\limits_a^{b_2}  \Big[ \big(\alpha \coth(\alpha s)
-\partial_1\big)h_k\big(s,\gamma_2(s)\big)\Big]
\Big[\big(\coth (\alpha s)
-\alpha^{-1}\partial _1\big)h_{\ell}\big(s,\gamma_2(s)\big)\Big] {\mathrm d}s \\
&+\!\!\!\!\!\int\limits_{\gamma_2(a)}^{\gamma_2(b_2)}\!\!\!\!
\coth\big(\alpha\gamma_2^{-1}(t)\big) 
\Big[\beta h_k\big(\gamma_2^{-1}(t),t\big)
h_{\ell}\big(\gamma_2^{-1}(t),t\big)+\beta^{-1}\partial_2
h_k\big(\gamma_2^{-1}(t),t\big)  
\partial_2 h_{\ell}\big(\gamma_2^{-1}(t),t\big)\Big]{\mathrm d}t, \\
A_{k,\ell}^{(4)}:=&  \frac 12  h_k\big(c,\gamma_0(c)\big) 
h_{\ell}\big(c,\gamma_0(c)\big)-\frac 12  h_k\big(b_2,\gamma_0(b_2)\big) 
h_{\ell}\big(b_2,\gamma_0(b_2)\big)\\
&+\int\limits_{b_2}^c \Big[ \alpha h_k\big(s,\gamma_0(s)\big)
h_{\ell}\big(s,\gamma_0(s)\big)+\alpha^{-1}\partial_1 h_k\big(s,\gamma_0(s)\big)
\partial _1h_{\ell}\big(s,\gamma_0(s)\big)\Big] {\mathrm d}s \\
&+\int\limits_{\gamma_0(c)}^{\gamma_0(b_2)}
\Big[\big(\beta+\partial_2\big)h_k\big(\gamma_0^{-1}(t),t\big)
\Big]\Big[\big(1+\beta^{-1}
\partial_2\big)h_{\ell}\big(\gamma_0^{-1}(t),t\big)\Big] {\mathrm d}t, \\
A_{k,\ell}^{(5)}:=&\iint\limits_G \big [\alpha\beta
h_k(s,t)h_{\ell}(s,t)+\alpha^{-1}\beta \partial_1h_k(s,t)\partial_1h_{\ell}(s,t)
+\alpha\beta^{-1} \partial_2h_k(s,t)\partial_2h_{\ell}(s,t) \\
&\phantom{+\iint\limits_G=}
+\alpha^{-1}\beta^{-1} \partial_1\partial_2h_k(s,t)
\partial_1\partial_2h_{\ell}(s,t) \big]\,{\mathrm d}s\, {\mathrm d}t,
\end{align*}
which equals \eqref{eq:eq2.10}.

At the end, using again \eqref{eq:eq4.5} and Remark \ref{integral} we have
\begin{equation*}
\zeta_k=\zeta_k^{(1)}+\zeta_k^{(2)}+\zeta_k^{(3)}+\zeta_k^{(4)}+\zeta_k^{(5)}
+\zeta_k^{(6)},  
\end{equation*}
where
\begin{align*}
\zeta_k^{(1)}:=&\,\frac{4{\mathrm e}^{2\alpha
    b_1+2\beta\gamma_{1,2}(b_1)}}{({\mathrm e}^{2\alpha b_1}-1)
  ({\mathrm e}^{2\beta\gamma_{1,2}(b_1)}-1)} 
  h_k\big(b_1,\gamma_{1,2}(b_1)\big) Z\big(b_1,\gamma_{1,2}(b_1)\big)\\
\zeta_k^{(2)}:=&\, \qlim _{\varrho\to 0} \frac 1{\varrho}
\int\limits_{b_1}^c \frac{4{\mathrm e}^{\alpha
    s+2\beta\gamma_1(s)}}{{\mathrm e}^{2\beta\gamma_1(s)}-1} 
(\alpha+\partial_1)h_k\big(s,\gamma_1(s)\big) \\
&\phantom{\, \qlim _{\varrho\to 0} \frac 1{\varrho}
\int\limits_{b_1}^c }
\times\bigg[\big({\mathrm e}^{2\alpha
  s}+\varrho\big)^{1/2}Z\Big(\frac{\log({\mathrm 
  e}^{2\alpha s}+\varrho)}{2\alpha},\gamma_1(s)\Big) -{\mathrm
  e}^{\alpha s}Z\big(s,\gamma_1(s)\big)\bigg]\,{\mathrm d}s,\\
\zeta_k^{(3)}:=&\, \int\limits_a^{b_1} \frac{4{\mathrm e}^{2\alpha
    s+2\beta\gamma_{1,2}(s)}}{({\mathrm e}^{2\alpha s}-1)
  ({\mathrm e}^{2\beta\gamma_{1,2}(s)}-1)} \big(\alpha \coth(\alpha s)
-\partial_1\big)h_k\big(s,\gamma_{1,2}(s)\big)Z\big(s,\gamma_{1,2}(s)\big)
{\mathrm d}s \\
&-\qlim _{\varrho\to 0} \frac 1{\varrho}
\int\limits_a^{b_1} \frac{4{\mathrm e}^{\alpha
    s+2\beta\gamma_{1,2}(s)}}{{\mathrm e}^{2\beta\gamma_{1,2}(s)}-1}
\big(\alpha \coth(\alpha s)-\partial_1)h_k\big(s,\gamma_{1,2}(s)\big) \\
&\phantom{\, \qlim _{\varrho\to 0} \frac 1{\varrho}
\int\limits_a^{b_1} }
\times \bigg[\big({\mathrm e}^{2\alpha
  s}+\varrho\big)^{1/2}Z\Big(\frac{\log({\mathrm 
  e}^{2\alpha s}+\varrho)}{2\alpha},\gamma_{1,2}(s)\Big) -{\mathrm
  e}^{\alpha s}Z\big(s,\gamma_{1,2}(s)\big)\bigg]\,{\mathrm d}s, \\
\zeta_k^{(4)}:=&\, \qlim _{\delta\to 0} \frac 1{\delta}
\int\limits_{\gamma_2(a)}^{\gamma_2(b_2)} 
\frac{4{\mathrm e}^{\beta t+2\alpha\gamma_2^{-1}(t)}}{{\mathrm
    e}^{2\alpha\gamma_2^{-1}(t)}-1} 
(\beta+\partial_2)h_k\big(\gamma_2^{-1}(t),t\big) \\
&\phantom{\, \qlim _{\varrho\to 0} \frac 1{\varrho}
\int\limits_{\gamma_2(a)}^{\gamma_2(b_2)} }
\times\bigg[\big({\mathrm e}^{2\beta
  t}+\delta\big)^{1/2}Z\Big(\gamma_2^{-1}(t),\frac{\log({\mathrm 
  e}^{2\beta t}+\delta)}{2\beta}\Big) -{\mathrm
  e}^{\beta t}Z\big(\gamma_2^{-1}(t),t\big)\bigg]\,{\mathrm d}t, \\
\zeta_k^{(5)}:=&\, \qlim _{\delta\to 0} \frac 1{\delta}
\int\limits_{\gamma_{1,2}(b_1)}^{\gamma_{1,2}(a)} 
\frac{4{\mathrm e}^{\beta t+2\alpha\gamma_{1,2}^{-1}(t)}}{{\mathrm
    e}^{2\alpha\gamma_{1,2}^{-1}(t)}-1} 
(\beta+\partial_2)h_k\big(\gamma_{1,2}^{-1}(t),t\big) \\
&\phantom{\, \qlim _{\varrho\to 0} \frac 1{\varrho}
\int\limits_{\gamma_{1,2}(b_1)}^{\gamma_{1,2}(1)} }
\times\bigg[\big({\mathrm e}^{2\beta
  t}+\delta\big)^{1/2}Z\Big(\gamma_{1,2}^{-1}(t),\frac{\log({\mathrm 
  e}^{2\beta t}+\delta)}{2\beta}\Big) -{\mathrm
  e}^{\beta t}Z\big(\gamma_{1,2}^{-1}(t),t\big)\bigg]\,{\mathrm d}t,
\end{align*}
\begin{align*}
\zeta_k^{(6)}:=&\,\qlim _{\varrho,\delta\to 0} \frac 4{\varrho \delta}
\iint\limits _G  {\mathrm e}^{\alpha s+\beta t}
(\alpha+\partial_1)(\beta+\partial_2)h_k\big(s,t\big)\bigg[{\mathrm e}^{\alpha
  s+\beta t}Z(s,t) \\
&\phantom{,\qlim _{\varrho,\delta\to 0} \frac 4{\varrho \delta}
\iint\limits _G }
-\!\big({\mathrm e}^{2\alpha
  s}\!+\!\varrho\big)^{1/2}{\mathrm e}^{\beta
  t}Z\Big(\frac{\log({\mathrm 
  e}^{2\alpha s}\!+\!\varrho)}{2\alpha},t\Big)\!-\!{\mathrm e}^{\alpha
  s}\big({\mathrm e}^{2\beta
  t}\!+\!\delta\big)^{1/2}Z\Big(s,\frac{\log({\mathrm 
  e}^{2\beta t}\!+\!\delta)}{2\beta}\Big)\\
&\phantom{,\qlim _{\varrho,\delta\to 0} \frac 4{\varrho \delta}
\iint\limits _G }
+\!\big({\mathrm e}^{2\alpha
  s}+\varrho\big)^{1/2}\big({\mathrm e}^{2\beta
  t}+\delta\big)^{1/2}Z\Big(\frac{\log({\mathrm 
  e}^{2\alpha s}+\varrho)}{2\alpha},\frac{\log({\mathrm 
  e}^{2\beta t}+\delta)}{2\beta}\Big)\bigg]{\mathrm d}s\,{\mathrm d}t,
\end{align*}
implying \eqref{eq:eq2.11}. \proofend

\bigskip
\noindent
{\bf Acknowledgments.} \ Research has been supported by 
the Hungarian  Scientific Research Fund under Grant No. OTKA
T079128/2009 and partially supported by T\'AMOP
4.2.1./B-09/1/KONV-2010-0007/IK/IT project. 
The project is implemented through the New Hungary Development Plan
co-financed by the European Social Fund, and the European Regional
Development Fund.

\end{document}